\documentclass[11pt]{article}
\usepackage[latin1]{inputenc}
\usepackage[T1]{fontenc}
\usepackage{amsmath,amssymb,amsthm}

\DeclareMathOperator{\coker}{coker}
 \DeclareMathOperator{\HH}{H}
 \DeclareMathOperator{\Spec}{Spec}
\DeclareMathOperator{\Hom}{Hom}
 \DeclareMathOperator{\Ext}{Ext}
\DeclareMathOperator{\diam}{diam}

\DeclareMathOperator{\ext}{ext}
\DeclareMathOperator{\Pic}{Pic}

\newtheorem{theorem}{Theorem} \newtheorem{lemma}[theorem]{Lemma}

\newtheorem{definition}[theorem]{Definition}
\newtheorem{example}[theorem]{Example}
\newtheorem{proposition}[theorem]{Proposition}

\newtheorem{remark}[theorem]{Remark}
\newtheorem{conjecture}[theorem]{Conjecture}

\usepackage{latexsym}            
\usepackage{amssymb}

\newcommand{\pp}{{\mathbb P}}

\newcommand\sE{{\mathcal E}}
\newcommand\sF{{\mathcal F}}
\newcommand\sG{{\mathcal G}}
\newcommand\sH{{\mathcal H}}
\newcommand\sI{{\mathcal I}}

\newcommand\sL{{\mathcal L}}

\newcommand\sN{{\mathcal N}}
\newcommand\sO{{\mathcal O}}

\newcommand{\proj}[1]
{ \mathchoice
           { {\mathbb P}^{#1} }
           { {\mathbb P}^{#1} }
           { {\mathbb P}^{#1} }
           { {\mathbb P}^{#1} }
         }

\begin{document}

\title{Moduli spaces of reflexive sheaves of rank 2}

\author{Jan O. Kleppe}


\date{\ }

\maketitle

\vspace*{-0.50in}
\begin{abstract}
\noindent 
Let $\sF$ be a coherent rank $2$ sheaf on a scheme $Y \subset \proj{n}$ of
dimension at least two. In this paper we study the relationship between the
functor which deforms a pair $(\sF,\sigma)$, $\sigma \in H^0(\sF)$, and the
functor which deforms the corresponding pair $(X,\xi)$ given as in the Serre
correspondence. We prove that the scheme structure of e.g. the moduli scheme
${\rm M_Y}(P)$ of stable sheaves on a threefold $Y$ at
$(\sF)$, 
and the scheme structure at $(X)$ of the Hilbert
scheme 
of curves on $Y$ are closely related. Using this relationship we get criteria
for the dimension and smoothness of $ {\rm M_{Y}}(P)$ at $(\sF)$, without
assuming $ {\Ext^2}(\sF ,\sF ) = 0$. For reflexive sheaves on $Y=\proj{3}$
whose deficiency module $M = H_{*}^1(\sF)$ satisfies $ {_{0}\!\Ext^2}(M ,M ) =
0 $ (e.g. of diameter at most 2), 
we get necessary and sufficient conditions of unobstructedness which coincide
in the diameter one case. The conditions are further equivalent to the
vanishing of certain graded Betti numbers of the free graded minimal
resolution of $H_{*}^0(\sF)$. It follows that every irreducible component of
${\rm M}_{\proj{3}}(P)$ containing a reflexive sheaf of diameter one is
reduced (generically smooth). We also determine a good lower bound for the
dimension of any component of ${\rm M}_{\proj{3}}(P)$ which contains a
reflexive stable sheaf with ``small'' deficiency module $M$.



 \end{abstract}
 \vspace*{-0.25in}
\thispagestyle{empty}

\section{Introduction and Main Results}
Let $Y \subset \proj{n}$ be an equidimensional, locally Cohen-Macaulay (CM),
closed subscheme of dimension at least two and let $\sF$ be a coherent rank
$2$ sheaf on $Y$. Let ${\rm Hilb}_{X/Y}$ be the local Hilbert functor of flat
deformations $X_S \subset Y \times S$, $S$ a local artinian $k$-algebra, of a
codimension 2 subscheme $X$ of $ Y$. An effective method of studying the
Hilbert scheme, ${\rm Hilb}^{p}(Y)$, of subschemes of $Y$ with Hilbert
polynomial $p$ with respect to smoothness, dimension and irreducibility at
$(X)$, is to look at other local deformation functors $\bf D$ over ${\rm
  Hilb}_{X/Y}$, ${\bf D} \rightarrow {\rm Hilb}_{X/Y}$, which allow a
surjective tangent map $t_{\bf D} \rightarrow t_{{\rm Hilb}_{X/Y}} =
H^0({\sN}_{X/Y})$, ${\sN}_{X/Y}=(\sI_{X/Y}/\sI_{X/Y}^2)^*$, and a
corresponding injective map of obstruction spaces. We consider such
deformation functors $\bf D$ which determine ${\rm Hilb}^{p}(Y)$ locally under
various assumptions. In particular we look to the functor of deformations of a
pair $(\sF,\sigma)$, or equivalently of the pair $(X,\xi)$ where $\xi$ is an
extension as in the Serre correspondence
\begin{equation} \label{serre} \xi \ \ ; \ \ 0 \rightarrow \sO_{Y}
  \stackrel{\sigma}{\longrightarrow} \sF \rightarrow \sI_{X/Y} \otimes \sL
  \rightarrow 0 \ . 
\end{equation} 
see \cite{H3}, \cite{HL}, \cite{Val}, \cite{Ver} and \cite{Vog} for the
existence of such extensions. Let ${\rm Def}_{\sF}$ (resp. ${\rm
  Def}_{\sF,\sigma}$) be the local deformation functor of flat deformations
$\sF_S$ of $\sF$ (resp. $\sO_{Y \times S} \stackrel{\sigma_S}{\longrightarrow}
\sF_S$ of $\sO_{Y} \stackrel{\sigma}{\longrightarrow} \sF$). Note that we have
an obvious forgetful map $p: {\rm Def}_{\sF,\sigma} \rightarrow {\rm
  Def}_{\sF}$. There exists also a natural projection $ {\rm Def}_{\sF,\sigma}
\rightarrow {\rm Hilb}_{X/Y}$ given by
$$(\sO_{Y \times S} \stackrel{\sigma_S}{\longrightarrow} \sF_S) \rightarrow
((\coker \sigma_S) \otimes_{\sO_{Y \times S}}(\sO_{Y \times S} \otimes_{\sO_Y}
{\sL}^{-1}))$$ which one may think of as given by the (relative) Serre
correspondence and the forgetful map $(X_S,\xi_S) \rightarrow (X_S)$. We
shortly write $\coker \sigma_S \otimes {\sL}^{-1}$ for $(\coker \sigma_S)
\otimes_{\sO_{Y \times S}}(\sO_{Y \times S} \otimes_{\sO_Y} {\sL}^{-1})$.
Using small letters for the dimensions, e.g. $ \ext^i(\sF, \sF) = \dim
\Ext^i(\sF, \sF)$ and $h^0(\sF)= \dim H^0(\sF)$ we prove
\begin{theorem} \label{localserre} With the notations above, suppose $Y$ is
  locally CM and equidimensional of dimension $ \dim Y \ge 2$ and such that
  $H^0(\sO_Y) \simeq k$ and $H^i(\sO_Y)=0$ for $i = 1$ and $2$. Moreover we
  suppose there exists an exact sequence \eqref{serre} where $X$ is
  equidimensional and locally CM of codimension $2$ in $Y$ and $\sI_{X/Y} =
  \ker (\sO_Y \rightarrow \sO_X)$. Then $ \Ext^1(\sI_{X/Y} \otimes \sL, \sF)$
  is the tangent space of $ {\rm Def}_{\sF,\sigma}$ and $\Ext^2(\sI_{X/Y}
  \otimes \sL, \sF)$ contains the obstructions of deforming $(\sF,\sigma)$.
  Moreover
  \\[-3mm]

  $\ \ \ \ (i) \ \ \ \ \ p: {\rm Def}_{\sF,\sigma} \rightarrow {\rm
    Def}_{\sF}$ is smooth (i.e. formally smooth) provided $H^1(\sF)=0, \ {
    and}$ \\[-4mm] 
  
  $\ \ \ \ (ii) \ \ \ \ q: {\rm Def}_{\sF,\sigma} \rightarrow {\rm
    Hilb}_{X/Y}$ is smooth provided $\Ext_{\sO_Y}^2(\sF, \sO_Y)=0 . $ \\
  Furthermore suppose $ H^1(\sF)=0$, $\Ext_{\sO_Y}^2(\sF, \sO_Y)=0 $, 
and that $\omega_Y$ is invertible. Then
  $$ \ext^1(\sF, \sF) - hom(\sF, \sF)+h^0(\sF) = h^0(\sN_{X/Y}) - 1 +
  h^0(\omega_X \otimes \omega_Y^{-1} \otimes {\sL}^{-1})- \sum_{i=0}^2 h^i(
  {\sL}^{-1}).$$ Suppose in addition that $\sF$ is stable (GM-stable) and
  $H^i( {\sL}^{-1})=0$ for $i=0,1,2$. Then
 $$  {\dim_{(\sF)}{\rm M_{Y}}}(P)+h^0(\sF) =
 \dim_{(X)}{\rm Hilb}^{p}(Y) + h^0(\omega_X \otimes \omega_Y^{-1} \otimes
 {\sL}^{-1})\ , \ \ {\rm and} \ $$ ${\rm M_{Y}}(P)$, the moduli scheme of
 stable sheaves with Hilbert polynomial $P$ on $Y$, is smooth at $(\sF)$ if
 and only if ${\rm Hilb}^{p}(Y)$ is smooth at $(X)$. Moreover $\sF$ is a
 generic sheaf of some component of ${\rm M_Y}(P)$ if and only if $X$ is
 generic in a corresponding component of ${\rm Hilb}^{p}(Y)$.
\end{theorem}

Let $F = H_{*}^0(\sF):= \oplus H^0(\sF(v))$, $M = H_{*}^1(\sF)$ and $E =
H_{*}^2(\sF)$. If $ {_0\!\Hom}(F , M) = 0$ and $Y$ is arithmetically
Cohen-Macaulay (ACM), then we show that the local graded deformation functors
of $F$ and of $(F, H_{*}^0(\sigma))$ are isomorphic to $ {\rm Def}_{\sF}$
and ${\rm Def}_{\sF,\sigma}$ respectively. We get the following variation of
Theorem~\ref{localserre}(i); that $p$ is smooth provided $ {_0\!\Hom}(F, M)
= 0$ and $Y$ is ACM. 

One 
may interpret the morphisms $p$ and $q$ in Theorem~\ref{localserre} as
corresponding to natural projections in an incidence correspondence of schemes
of corepresentable functors, connecting ${\rm M_Y}(P)$ closely to ${\rm
  Hilb}^{p}(Y)$. 
Under the assumptions of Theorem~\ref{localserre} the projections are smooth
of known fiber dimension. Since the fiber dimensions are easy to see and the
Serre correspondence is well understood (\cite{H3}), related arguments as in
the theorem are used in the literature, especially to compute dimensions of or
describe very specific moduli schemes (e.g. \cite{C2}, \cite{H2}, \cite{HS},
\cite{IM}, \cite{Mir3}, \cite{St}, \cite{Gu} and see \cite{Ver}, sect.\! 4 for
results and a discussion). It is, however, under the mere assumptions of {\rm
  (i)} and {\rm (ii)} above we are able to precisely see that the scheme
structures of ${\rm M_Y}(P)$ and ${\rm Hilb}^{p}(Y)$ are ``the same''. To
apply Theorem~\ref{localserre} we neither need $ H^1(\sN_{X/Y})=0$, nor
$\Ext_{\sO_Y}^2(\sF, \sF)=0 $ to prove the smoothness of the moduli schemes.
This, we think, significantly distinguishes our theorem from the results and
the proofs of 
the mentioned papers. For the complete picture we have no better reference
that an old preprint of the author (\cite{K2}, for the case $Y=\proj{3}$) and
the paper \cite{GMir} which explicitly makes use of (without proofs) and
slightly extends the results of \cite{K2}, and we therefore include full
proofs.

As an application we prove several results concerning smoothness and dimension
of the moduli space, ${\rm M}_{\proj{3}}(c_1,c_2,c_3)$, of stable {\it
  reflexive} sheaves of rank 2 with Chern classes $c_1,c_2$ and $c_3$ on
$\proj{3}$. In some cases, especially for $c_3=0$ or small $c_2$ or large
$c_3$, one knows the answer, e.g. see \cite{AO}, \cite{C1}, \cite{C2},
\cite{ES}, \cite{GMir}, \cite{H2}, \cite{H3}, \cite{HS}, \cite{Mir2},
\cite{Mir3} and \cite{MSS}. Much is still unknown about ${\rm
  M}_{\proj{3}}(c_1,c_2,c_3)$, see \cite{RVV} for an overview of recent
research. Let
$${\rm ed}(\sF)= \ext_{\sO_Y}^1(\sF, \sF) - \ext_{\sO_Y}^2(\sF, \sF)$$ (using
small letters for the dimension of the global Ext-group). If $\sF$ is stable,
then ${\rm ed}(\sF)$ is sometimes called the ``expected dimension'' of
${\rm M_Y}(P)$ at $(\sF)$ and ${\rm ed}(\sF)=8c_2-2c_1^2-3$ if $Y=\proj{3}$. 
We prove that ${\rm M}_{\proj{3}}(c_1,c_2,c_3)$ is smooth at $(\sF)$, i.e.
that $\sF$ is unobstructed, and we find $\dim_{(\sF)}{\rm
  M}_{\proj{3}}(c_1,c_2,c_3) $ provided we have sufficient vanishing of
${_v\!\Hom_R}(F,M)$ and ${_v\!\Hom_R}(M , E)$ for $v = 0$ and $ -4$
(Theorem~\ref{introrefl1}). This result generalizes \cite{Mir2}, which gives
the complete answer for $M=0$. 
Let  ${_0\!\Ext_R^2}(M , M) = 0$.  Using that the composition
\begin{equation} \label{cup} \eta \ : \ {_0\!\Hom_R}(F,M) \times
  {_0\!\Hom_R}(M , E) \longrightarrow {_0\!\Hom_R}(F , E) \ ,
\end{equation}
commutes with the cup product, we show that $\sF$ is obstructed if $\eta \ne
0$ (cf. \cite{F}, \cite{krao}, \cite{W1}). 
Thanks to this result we get that the sufficient conditions of
unobstructedness of Theorem~\ref{introrefl1} are close (resp. equivalent) to
being {\it necessary} conditions provided the diameter of $M$ is small (resp.
one). Since we can substitute the non-vanishing of the $\Hom$-groups of
Theorem~\ref{introrefl1} by the non-triviality of certain graded Betti numbers
appearing in the minimal resolution,
$$ 0 \rightarrow \bigoplus_i \sO_{\proj{}}(-i)^{\beta_{3,i}} \rightarrow
\bigoplus_i \sO_{\proj{}}(-i)^{\beta_{2,i}} \rightarrow \bigoplus_i
\sO_{\proj{}}(-i)^{\beta_{1,i}} \rightarrow \sF \rightarrow 0 \ ,
 $$ of $\sF$ and we 
 explicitly compute ${_0\!\hom_R}(F,E)$ (see Remark~\ref{explici}), we get

\begin{theorem} \label{introrefl3} Let $\sF$ be reflexive of rank $2$ on
  $\proj{3}$ and suppose $M \neq 0$ is $r$-dimensional of diameter $1$ and
  concentrated in degree $c$. 
  Then $\sF$ is obstructed if and only if
  $$
  \beta_{1,c} \cdot \beta_{2,c+4} \neq 0 \ \ \ {\rm or} \ \ \ \beta_{1,c+4}
  \cdot \beta_{2,c+4} \neq 0 \ \ \ {\rm or} \ \ \ \beta_{1,c} \cdot
  \beta_{2,c}\neq 0 \ . $$ Moreover if $\sF$ is an unobstructed stable sheaf,
  then the dimension of the moduli scheme ${\rm M}_{\proj{3}}(c_1,c_2,c_3)$ at
  $(\sF)$ is
\begin{equation*}
  \dim_{(\sF)}{\rm M}_{\proj{3}}(c_1,c_2,c_3) =  {\rm ed}(\sF)+
  {_0\!\hom_R}(F,E) + 
  r(\beta_{1,c+4} +  \beta_{2,c}). 
\end{equation*}
\end{theorem}

It follows that every irreducible component of ${\rm
  M}_{\proj{3}}(c_1,c_2,c_3)$ whose generic sheaf $\sF$ satisfies $\diam M =
1$ is {\it reduced} (i.e. generically smooth). If $\diam M = m$ we give
examples of moduli spaces ${\rm M}_{\proj{3}}(c_1,c_2,c_3)$ containing a
non-reduced component for every integer $m \ge 3$.
If $\diam M = 2$ we conjecture that the corresponding component of ${\rm
  M}_{\proj{3}}(c_1,c_2,c_3)$ is generically smooth. We also give a new
formula for the dimension of {\it any} generically smooth irreducible
component of ${\rm M}_{\proj{3}}(c_1,c_2,c_3)$ (Theorem~\ref{thmeulerIF}).
Even though some of the results of this paper may have a direct proof in which
the condition ``reflexive'' is replaced by ``torsionfree'', we have chosen
just to use Theorem~\ref{localserre} and the corresponding results for ${\rm
  Hilb}^{p}(\proj{3})$.

The main results of this paper were lectured at the workshop on "Vector
Bundles and Low Codimensional Subvarieties'' at Trento, in September 2006. The
author thanks the organizers for their hospitality. Theorem~\ref{localserre}
for $Y=\proj{3}$, Example~\ref{can} and Example~\ref{mum} were the main
results (except for the dimension statements) in an old preprint of the author
(\cite{K2}, available at \text{http://www.iu.hio.no/\textasciitilde
  jank/papers.htm}). 
Moreover I heartily thank prof. O. A.
Laudal at Oslo, prof. S.A. Strømme at Bergen and prof. R.M. Mir\'o-Roig at
Barcelona for interesting discussions and comments.

\subsection {Notations and terminology} $R = k[X_0,X_1,...,X_n]$ is a graded
polynomial ring over an algebraically closed $k$ of arbitrary characteristic
with the standard grading, $ \mathfrak m = (X_0,..,X_n)$ and $Y \subset
\proj{n}$ is a closed equidimensional, locally Cohen-Macaulay (CM) subscheme.
We keep the other notations of the introduction. A {\it curve} $X$ in
$\proj{n}$ (resp. in $Y $) is an {\it equidimensional, locally CM} subscheme
of $\proj{}:=\proj{n}$ (resp. of $Y$) of dimension one with sheaf ideal
$\sI_X$ (resp. $\sI_{X/Y}$) and normal sheaf $\sN_X =
{\sH}om_{\sO_{\proj{}}}(\sI_X,\sO_X)$ (resp. $\sN_{X/Y} =
\Hom_{\sO_{Y}}(\sI_{X/Y},\sO_X)$ in $Y$). $X$ is {\it unobstructed\/} if the
Hilbert scheme is smooth at the corresponding point $(X)= (X \subseteq
\proj{n})$, otherwise $X$ is {\it obstructed}. The Hilbert scheme of space
curves of degree $d$ and arithmetic genus $g$ is denoted by $\HH(d,g)$. If
$\sF$ is a coherent $\sO_{Y}$-Module, we let $H^i(\sF) =
H^i(Y,\sF)$ 
and $h^i(\sF) = \dim H^i(\sF)$, and we denote by $\chi(\sF) = \Sigma (-1)^i
h^i(\sF)$. 
Then $I_X:= H_{*}^0(\proj{},\sI_X)$ is the saturated homogeneous ideal of $X$
in $\proj{n}$.

Let $M = M(\sF)$ be the deficiency module $H_{*}^1(\sF)$. $\sF$ is said to be
unobstructed if the hull of the local deformation functor, ${\rm Def}_{\sF}$,
is smooth. By stable we mean GM-stable, i.e. stable in the sense of Gieseker
and Maruyama in which the Hilbert polynomial (and not the $1^{st}$ Chern
class) is used to define stability (see \cite{HL}, ch. I). Thus a stable $\sF$
is unobstructed iff ${\rm M_Y}(P)$, the moduli scheme of stable sheaves with
Hilbert polynomial $P$ on $Y$, is smooth at $(\sF)$ (\cite{HL}, Thm. 4.5.1).
The two concepts of stability are the same if $Y=\proj{3}$ (\cite{H3}, Rem.\!
3.1.1). Stable sheaves are simple, i.e. $\Hom(\sF,\sF) \simeq k$ (\cite{HL},
Cor.\! 1.2.8). A coherent sheaf $\sF$ is reflexive iff $\sF \simeq \sF^{**}$
where $\sF^{*}={\sH}om(\sF,-)$ (see \cite{H3}). In the case $Y$ is a smooth
threefold, we denote by ${\rm M_{Y}}(c_1,c_2,c_3)$ the moduli scheme of stable
{\it reflexive} sheaves of rank 2 on $Y$ with Chern classes $c_1,c_2$ and
$c_3$. Thus ${\rm M_{Y}}(c_1,c_2,c_3)$ is open in ${\rm M_{Y}}(P)$. For any
$\sF$ of ${\rm M_{Y}}(c_1,c_2,c_3)$ there exists an exact sequence
\eqref{serre} after replacing $\sF$ by some $\sF(t)$, due to a theorem of
Hartshorne. In fact \eqref{serre} defines a correspondence between the sets of
pairs $(\sF,s)$ and pairs $(X,\xi)$ provided $\xi$ generates $\omega_X \otimes
\omega_Y^{-1} \otimes {\sL}^{-1}$ except at finitely many points (the
Hartshorne-Serre correspondence, see \cite{H3}, Thm.\! 4.1 and \cite{Val},
Thm.\! 1).

A sheaf $\sF$ of rank $2$ on $\proj{3}$ is said to be Buchsbaum if $\mathfrak
m \cdot M(\sF) = 0$. We define the diameter of $M(\sF)$ (or of $\sF$) by
$\diam M(\sF) = c- b+1 $ where $b=\min \{n \arrowvert h^1(\sF(n)) \neq 0 \}$
and $c=\max \{n \arrowvert h^1(\sF(n)) \neq 0 \}$ and by $\diam M(\sF) = 0$ if
$M(\sF)=0$. The diameter of a curve $C$, $\diam M(C)$, is correspondingly
defined. 
A curve in a sufficiently small open irreducible subset of $\HH(d,g)$ (small
enough to satisfy all the openness properties which we want to pose) is called
a {\it generic} curve of $\HH(d,g)$, and accordingly, if we state that a
generic curve has a certain property, then there is a non-empty open
irreducible 
subset of $\HH(d,g)$ of curves having this property. A {\it generization} $C'
\subseteq \proj{3}$ of $C \subseteq \proj{3}$ in 
$\HH(d,g)$ is a the generic curve of some irreducible subset of $\HH(d,g)$
containing $(C)$. In the same way we use the word generic and generization for
a stable sheaf. By an irreducible component of $\HH(d,g)$ or ${\rm
  M}_{\proj{3}}(c_1,c_2,c_3)$ we always mean a {\it non-embedded} irreducible
component.

For any graded R-module $N$ we have the right derived functors $H_{\mathfrak
  m}^i(N)$ and ${_v\!\Ext_{\mathfrak m}^i}(N , -)$ of $\Gamma_{\mathfrak m}(N)
= \sum_v \ker(N_v \rightarrow \Gamma( \proj{}, \tilde N (v)))$ and
$\Gamma_{\mathfrak m}(\Hom_R(N , -))_v$ respectively (cf. \cite{SGA2}, exp.
VI). 
We use small letters for the $k$-dimension and subscript $v$ for the
homogeneous part of degree $v$, e.g. ${_v\!\ext_{\mathfrak m}^i}(N_1,N_2) =
\dim {_v\!\Ext_{\mathfrak m}^i}(N_1,N_2)$, for graded $R$-modules $N_i$ of
finite type. There is a spectral sequence (\cite{SGA2}, exp. VI)
\begin{equation}\label{mspect}
E_2^{p,q} = {_v\!\Ext_R^p }(N_1 , H_{\mathfrak m}^q (N_2 )) \Rightarrow
{_v\!\Ext_{\mathfrak m}^{p+q} }(N_1, N_2 )
\end{equation}
($\Rightarrow$ means ``converging to'') and a duality isomorphism (\cite{K5},
Thm. 1.1);
\begin{equation}\label{mduality}
{_v\!\Ext^{i+1}_{\mathfrak m}} (N_2 ,N_1 ) \simeq \hbox{$_{-v}
  \Ext_R^{n-i} (N_1 ,N_2(-n-1) )^\vee$}
\end{equation}
where $(-)^\vee = \Hom_k( - , k)$, generalizing the Gorenstein duality
$_{v}\!H_{\mathfrak{m}}^{i+1}(M) \simeq \
_{-v}\!\Ext_R^{n-i}(M,R(-4))^{\vee}.$ These groups fit into a long exact
sequence (\cite{SGA2}, exp.\ VI)
\begin{equation}\label{seq}
\rightarrow {_v\!\Ext_{\mathfrak m}^i }(N_1 ,N_2 ) \rightarrow {_v\!\Ext_R^i }
(N_1 ,N_2 ) \rightarrow \Ext_{{\mathcal O}_{\proj{}}}^i (\tilde N_1
,\tilde N_2 (v)) \rightarrow  {_v\!\Ext_{\mathfrak m}^{i+1} }(N_1 ,N_2 )
\rightarrow
\end{equation} 
which e.g. relates the deformation theory of $X \subseteq \proj{3}$, described
by $H^{i-1}({\sN_X}) \simeq \Ext_{{\mathcal O}_{\proj{}}}^i (\sI_X,\sI_X)$
for $i = 1,2$, to the deformation theory of the homogeneous ideal $I = I_X$
(or equivalently of $A=R/I$), described by ${_0\!\Ext_R^i } (I_X,I_X)$, in the
following exact sequence
\begin{equation}\label{mseq}
{_v\!\Ext_R^1 }(I,I ) \hookrightarrow  H^0({\sN_C}(v)) \rightarrow
{_v\!\Ext_{\mathfrak m}^2}(I,I ) \stackrel{\alpha}{\longrightarrow}
{_v\!\Ext_R^2 }(I,I ) \rightarrow 
H^1({\sN_C}(v)) \rightarrow {_v\!\Ext_{\mathfrak m}^3 }(I ,I ) \rightarrow 0
\ 
\end{equation}
(cf. \cite{krao}, section 2). Let $M(X)= H_{\mathfrak m}^2(I)$. Note that, in
this situation, C. Walter proved in \cite{W2}, Thm. 2.3 that the map $\alpha:
{_v\!\Ext_{\mathfrak m}^2 }(I ,I ) \simeq {_v\!\Hom_R}(I ,H_{\mathfrak m}^2(I)
) \rightarrow {_v\!\Ext_R^2 }(I,I )$ of \eqref{mseq} factorizes via $
{_v\!\Ext_R^2 }(M(X), M(X))$ in a natural way, the factorization is in fact
given by a certain edge homomorphism of the spectral sequence \eqref{mspect}
with $N_1 = M$, $N_2=I$ and $p+q = 4$ (cf. \cite{F}, Thm\! 2.5). We
frequently refer to \cite{krao} and all results we use from \cite{krao} are
true without the characteristic zero assumption of the field quoted for that
paper.
%

\section{The scheme structure in the Serre correspondence}
In this section we will prove the basic Theorem~\ref{localserre} and its
variations. Moreover we give some applications and examples of moduli schemes
${\rm M_Y}(c_1,c_2,c_3)$ in the case $Y=\proj{3}$. In particular we show that
some ${\rm M}_{\proj{3}}(c_1,c_2,c_3)$ contains a non-reduced component.

The local deformation functors ${\rm Def}_{\sF}$, ${\rm Def}_{\sF,\sigma}$ and
${\rm Hilb}_{X/Y}$ of the Introduction were defined on the category
$\underline{l}$ whose objects are local artinian $k$-algebras $S$ with residue
field $k$ and whose morphisms are homomorphisms of local rings over $k$.
Associated to \eqref{serre} there is also another local deformation functor on
$\underline{l}$:
\begin{equation*} {\rm  Def}_{X/Y,\xi}(S) = \{\left (X_S  \subset Y_S, \xi_S)\
    \arrowvert \ (X_S \subset Y_S) \in {\rm Hilb}_{X/Y}(S)\ {\rm and}\ \xi_S
    \otimes_S k = \xi \right\}  
  \end{equation*}
  where $Y_S:= Y \times S$ (we shortly write $ Y \times S$ for $ Y \times
  \Spec(S)$), $\sL_S:=\sO_{Y_S} \otimes_{\sO_Y} {\sL}$ and $ \xi_S \in
  \Ext^1(\sI_{X_S/Y_S} \otimes \sL_S, \sO_{Y_S})$. Now
  Theorem~\ref{localserre} follows easily from the following basic result
  about the relationship of dimensions and the scheme structures in the Serre
  correspondence.

  \begin{theorem} \label{localserree} Let $Y \subset {\proj{n}}$ be locally CM
    and equidimensional of dimension $ \dim Y \ge 2$ such that $H^0(\sO_Y)
    \simeq k$ and $H^i(\sO_Y)=0$ for $i = 1,2$. Moreover suppose there exists
    an exact sequence \eqref{serre} where $X$ is equidimensional and locally
    CM of codimension $2$ in $Y$ and $\sI_{X/Y} = \ker (\sO_Y \rightarrow
    \sO_X)$. Then  \\[-2mm]

    {\rm (a)} $ \Ext^1(\sI_{X/Y} \otimes \sL, \sF)$ is the tangent space of $
    {\rm Def}_{\sF,\sigma}$ and $\Ext^2(\sI_{X/Y} \otimes \sL, \sF)$ contains
    the obstructions of deforming $(\sF,\sigma)$. Moreover $ {\rm
      Def}_{\sF,\sigma} \simeq {\rm Def}_{X/Y,\xi}$ are isomorphic on
    $\underline{l}$ and
    \\[-3mm]

  $\ \ \ \ (i) \ \ \ \ \ p: {\rm Def}_{\sF,\sigma} \rightarrow {\rm
    Def}_{\sF}$ is smooth (i.e. formally smooth) provided $H^1(\sF)=0, \ {
    and}$ \\[-4mm] 
  
  $\ \ \ \ (ii) \ \ \ \ q: {\rm Def}_{\sF,\sigma} \rightarrow {\rm
    Hilb}_{X/Y}$ is smooth provided $\Ext_{\sO_Y}^2(\sF, \sO_Y)=0 . $ 
  \\[-1mm]

  {\rm (b)} Suppose $ H^1(\sF)=0$, $\Ext_{\sO_Y}^2(\sF, \sO_Y)=0 $,
and that $\omega_Y$ is invertible. Then
  $$ \ext^1(\sF, \sF) - hom(\sF, \sF)+h^0(\sF) = h^0(\sN_{X/Y}) - 1 +
  h^0(\omega_X \otimes \omega_Y^{-1} \otimes {\sL}^{-1})- \sum_{i=0}^2
  (-1)^{i} h^i( {\sL}^{-1}).$$ Suppose in addition that $\sF$ is stable and
  $H^i( {\sL}^{-1})=0$ for $i=0,1,2$. Then
 $$  {\dim_{(\sF)}{\rm M_{Y}}}(P)+h^0(\sF) =
 \dim_{(X)}{\rm Hilb}^{p}(Y) + h^0(\omega_X \otimes \omega_Y^{-1} \otimes
 {\sL}^{-1}) . \ $$ It follows that ${\rm M_{Y}}(P)$ is smooth at $(\sF)$ if
 and only if ${\rm Hilb}^{p}(Y)$ is smooth at $(X)$. Furthermore $\sF$ is a
 generic sheaf of some component of ${\rm M_Y}(P)$ if and only if $X$ is
 generic in a corresponding component of ${\rm Hilb}^{p}(Y)$.
\end{theorem}

\begin{remark} \label{remver2} Under the assumptions of
  Theorem~\ref{localserree} {\rm (a)}, we get that $H^1(\sF) \simeq
  H^1(\sI_{X/Y} \otimes \sL)$ and $\Ext_{\sO_Y}^2(\sF, \sO_Y) \simeq
  \Ext^2(\sI_{X/Y} \otimes \sL, \sO_Y)$ by using \eqref{serre}. Moreover, if
  the dualizing sheaf $\omega_Y$ is invertible (i.e. $Y$ locally Gorenstein),
  then $$\Ext_{\sO_Y}^2(\sF, \sO_Y)^{\vee} \simeq H^{\dim Y -2}(\sF \otimes
  \omega_Y) \ \ {\rm and} \ \ \Ext^2(\sI_{X/Y} \otimes \sL, \sO_Y)^{\vee}
  \simeq H^{\dim Y -2}(\sI_{X/Y} \otimes \sL \otimes \omega_Y).$$
\end{remark}

In the case $Y$ is locally Gorenstein and the immersion $Y \hookrightarrow
\proj{n}$ induces an isomorphism $ {\mathbb Z}=\Pic(\proj{n}) \simeq \Pic(Y)$,
we will use this isomorphism to look upon the first Chern class $c_1$ as an
integer, i.e. $\sL \simeq \wedge^2 \sF = \sO_Y(c_1)$. Moreover put $\omega_Y =
\sO_Y(e)$. 
By Remark~\ref{remver2} Theorem~\ref{localserree} immediately implies

\begin{theorem} \label{localserre1} Suppose $Y \subset \proj{n}$ is locally
  Gorenstein and equidimensional of dimension $ \dim Y \ge 2$ such that
  $H^0(\sO_Y) \simeq k$, $H^i(\sO_Y)=0$ for $i = 1,2$ and $
  \Pic(\proj{n}) \simeq \Pic(Y)$. Moreover we suppose there is an exact
  sequence of $\sO_Y$-Modules
  \begin{equation} \label{serre1} \xi \ \ ; \ \ 0 \rightarrow \sO_{Y}
    \stackrel{\sigma}{\longrightarrow} \sF \rightarrow \sI_{X/Y}(c_1)
    \rightarrow 0 \ 
\end{equation} 
where $X$ is equidimensional and locally CM of codimension $2$ in $Y$ and
$\sI_{X/Y} = \ker (\sO_Y \rightarrow \sO_X)$. Then $ \Ext^1(\sI_{X/Y}(c_1),
\sF)$ is the tangent space of $ {\rm Def}_{\sF,\sigma}$ and
$\Ext^2(\sI_{X/Y}(c_1), \sF)$ contains the obstructions of
deforming $(\sF,\sigma)$. Moreover
\\[-3mm]

$\ \ \ \ (i) \ \ \ \ \ p: {\rm Def}_{\sF,\sigma} \rightarrow {\rm Def}_{\sF}$
is smooth provided $H^1(\sI_{X/Y}(c_1))=0 \ , \ {
  and}$ \\[-4mm]
  
$\ \ \ \ (ii) \ \ \ \ q: {\rm Def}_{\sF,\sigma} \rightarrow {\rm Hilb}_{X/Y}$
is smooth provided $H^{\dim Y-2}(\sI_{X/Y}(c_1+e))=0 \ . $ \\ Furthermore if
$H^i(\sO_Y(-c_1))=0$ for $i=0,1,2$, $H^1(\sI_{X/Y}(c_1))=0$, $H^{\dim
  Y-2}(\sI_{X/Y}(c_1+e))=0$ and $\sF$ is a stable sheaf, then
  $$ \ext^1(\sF, \sF)+h^0(\sF) = h^0(\sN_{X/Y})+ h^0(\omega_X(-c_1-e))\ , \
  \ {\rm  and}$$
 $$  {\dim_{(\sF)}{\rm M_{Y}}}(P)+h^0(\sF) =
 \dim_{(X)}{\rm Hilb}^{p}(Y) + h^0(\omega_X(-c_1-e)) .$$ Finally $\sF$ is a
 generic sheaf of some component of ${\rm M_Y}(P)$ if and only if $X$ is
 generic in a corresponding component of ${\rm Hilb}^{p}(Y)$.
\end{theorem}

Somehow we may look upon the theorem as the Hartshorne-Serre correspondence
for flat families. We don't, however, need $\sF$ to be reflexive (only
torsionfree as one may easily deduce from \eqref{serre1}).

We shortly return to the proof of Theorem~\ref{localserree}. Firstly we give
an example to see that conditions as in Theorem~\ref{localserre1} are
necessary for comparing the structure of $\HH(d,g)$ and ${\rm
  M}_{\proj{3}}(c_1,c_2,c_3)$, while the same example ``twisted'' leads to a
non-reduced component of ${\rm M}_{\proj{3}}(c_1,c_2,c_3)$ once the conditions
of the theorem are
satisfied. 
Below we will need the following result (\cite{K3}, Prop.\! 3.2). Let $C$ and
$X$ be two space curves which are algebraically linked by a complete
intersection of two surfaces of degrees $f$ and $g$ (a c.i. of type $(f,g)$),
see \cite{Mig} for the theory on linkage. If we suppose
\begin{equation} \label{mli}
H^1(\sI_C(v)) = 0 \ \ \ {\rm for} \ v= f , g , f - 4 \ {\rm and} \  g - 4,
\end{equation}
then $C$ is unobstructed (resp. generic) iff $X$ is unobstructed (resp.
generic), and we have
$$\dim_{(C)}\HH(d,g)+ h^0(\sI_C(f)) + h^0(\sI_C(g))= \dim_{(X)}\HH(d',g')+
h^0(\sI_X(f)) + h^0(\sI_X(g)).$$
\begin{example} \label{can} The generic curve $C$ of Mumford's well-known
  example of a non-reduced component of $\HH(14,24)$ satisfies $H^1(\sI_C(v))
  = 0$ for $v \neq 3,4,5$ (\cite{Mu}). Moreover there is a c.i. of type $(6,6)$
  containing $C$ whose linked curve is smooth. Hence by the result mentioned
  in \eqref{mli} the linked curve is the general curve $X$ of a non-reduced
  component of $\HH(22,56)$ of dimension $88$. We leave to the reader to
  verify that $X$ is subcanonical ($\omega_X \simeq \sO_X(5)$) and satisfies
  $H^1(\sI_X(v)) = 0$ for $v \neq 3,4,5$.

  (a) If we take a general element of $H^0(\sO_X) \simeq H^0(\omega_X(-5))
  \simeq \Ext^1(\sI_{X},\sO_{\proj{3}}(-9))$, we get an extension
  \begin{equation} \label{serre2} \xi \ \ ; \ \ 0 \rightarrow \sO_{\proj{3}}
    \stackrel{\sigma}{\longrightarrow} \sE \rightarrow \sI_{X}(9) \rightarrow
    0 \
\end{equation} 
in which $\sE$ is a stable vector bundle with $c_1 = 9$ and $c_1(\sE(-5)) =
-1$, $c_2(\sE(-5)) = 2$. It is well known that ${\rm M}_{\proj{3}}(-1,2,0)$ is
smooth \cite{HS}, i.e. $\sE$ is unobstructed while $X$ is obstructed. The
assumption $H^1(\sI_{X}(c_1+e))=0$ of Theorem~\ref{localserre1} is, however,
not satisfied. Indeed $H^1(\sI_{X}(c_1+e))=H^1(\sI_{X}(5)) \ne 0 . $

(b) If we take a general global section of $\sO_X(3) \simeq \omega_X(-2)$ we
get an extension
\begin{equation} \label{serre3} \xi \ \ ; \ \ 0 \rightarrow \sO_{\proj{3}}
  \stackrel{\sigma}{\longrightarrow} \sF \rightarrow \sI_{X}(6) \rightarrow 0
  \
\end{equation}
in which $\sF$ is a stable reflexive sheaf belonging to ${\rm
  M}_{\proj{3}}(6,22,66) \simeq {\rm M}_{\proj{3}}(0,13,66)$. Since all
assumptions of Theorem~\ref{localserre1} are satisfied, we conclude that $\sF$
is the general point of a non-reduced component of ${\rm
  M}_{\proj{3}}(0,13,66)$ of dimension $-h^0(\sF)+\dim_{(X)}\HH(22,56)+
h^0(\omega_X(-2))=-8+88+21=101$. Note that we in this case have ${\rm
  ed}(\sF)=8c_2-2c_1^2-3=101$, i.e. the component is non-reduced of the least
possible dimension.
\end{example}
\begin{example} \label{mum} Here we apply Theorem~\ref{localserre1} directly
  to Mumford's example of a generic obstructed curve $C \in \HH(14,24)$.

(a) If we take a general element of
  $H^0(\omega_C(2)) \simeq \Ext^1(\sI_{C},\sO_{\proj{3}}(-2))$, we get an
  extension
  \begin{equation} \label{serre4} \xi \ \ ; \ \ 0 \rightarrow \sO_{\proj{3}}
    \stackrel{\sigma}{\longrightarrow} \sF \rightarrow \sI_{C}(2) \rightarrow
    0 \
\end{equation} 
in which $\sF$ is a stable reflexive sheaf with $c_1 = 2$ and $c_1(\sF(-1)) =
0$, $c_2(\sF(-1)) = 13$. The assumptions $H^1(\sI_{C}(c_1))=0$,
$H^{1}(\sI_{C}(c_1-4))=0$ of Theorem~\ref{localserre1} are satisfied and we
get a non-reduced component of ${\rm M}_{\proj{3}}(0,13,74)$ of dimension
$-h^0(\sF)+\dim_{(C)}\HH(14,24)+ h^0(\omega_C(2))=-1+56+51=106$.

(b) If we take a general global section of $H^0( \omega_C(3)) \simeq
\Ext^1(\sI_{C},\sO_{\proj{3}}(-1))$ we get by Theorem~\ref{localserre1} (with
$c_1=1$) an extension where $\sF$ is a stable reflexive sheaf belonging to a
non-reduced component of ${\rm M}_{\proj{3}}(1,14,88) \simeq {\rm
  M}_{\proj{3}}(-1,14,88)$ of dimension $-1+56+65=120$.

(c) If we take a general global section of $H^0( \omega_C(-2)) \simeq
\Ext^1(\sI_{C},\sO_{\proj{3}}(-6))$ we get by Theorem~\ref{localserre1} (with
$c_1=6$) an extension where $\sF$ is a {\rm semistable obstructed} reflexive
sheaf belonging to the moduli space of semistable sheaves ${\overline {\rm
    M}_{\proj{3}}}(6,14,18) \simeq {\overline {\rm M}_{\proj{3}}}(0,5,18)$.
Even though $\sF$ is obstructed, i.e. the hull of the local deformation
functor is singular, we don't yet know the hulls precise relationship to the
local ring $O_{{\overline M},(\sF)}$ of $ {\overline {\rm
    M}_{\proj{3}}}(0,5,18)$ at $(\sF)$ and we are not able to state whether
$O_{{\overline M},(\sF)}$ is singular or not.
\end{example}
\begin{proof}[Proof (of  Theorem~\ref{localserree})] {\rm (a)}
  Using Laudal's results (\cite{L1}) for the local deformation functor of
  deforming a category, we {\it claim} that $ \Ext^1(\sI_{X/Y} \otimes \sL,
  \sF)$ is the tangent space of ${\rm Def}_{\sF,\sigma}$ and that
  $\Ext^2(\sI_{X/Y} \otimes \sL, \sF)$ contains the obstructions of deforming
  $(\sF,\sigma)$. Indeed by \cite{L}, §2 there is a spectral sequence
  \begin{equation*}
    E_2^{p,q} = {\displaystyle\lim_{\longleftarrow}}^{(p)}
    \left\{
    \begin{array}{c}
      \Ext^q(\sF,\sF) \ \ \ \ \ \Ext^q(\sO_Y, \sO_Y)   \\
       \searrow^{\alpha^q} \ \ \ \ \ \ \
      \swarrow\ \  \ \ \ \ \  
      \\ \Ext^q(\sO_Y,\sF)  \ \ \ \ \ \ 
    \end{array}
    \right\}
  \end{equation*}
converging to an algebra cohomology group $A^{(.)}$ for which $A^1$ is the
tangent space of ${\rm Def}_{\sF,\sigma}$ and $A^2$ contains the obstructions
of deforming $(\sF,\sigma)$. Since $ E_2^{p,q}=0$ for $p \ge 2$, we get the
exact sequence
$$0 \rightarrow E_2^{1,q-1} \rightarrow A^{q} \rightarrow E_2^{0,q}
\rightarrow 0 \ .$$ Moreover $ \Ext^q(\sO_Y, \sO_Y)=0$ for $0 < q < 3$ by the
assumption $H^i(\sO_Y)=0$ for $i = 1, 2$, and we get $ E_2^{0,q}=\ker
\alpha^q$ and $E_2^{1,q}=\coker \alpha^q$ for $q >0$. Observe also that $
E_2^{1,0}= \coker \alpha^0$ because $k \simeq H^0(\sO_Y)\simeq \Hom(\sO_Y,
\sO_Y) \subseteq \Hom(\sF, \sF)$. We therefore have an exact sequence
$$ 0 \rightarrow  \coker \alpha^{q-1} \rightarrow A^{q} \rightarrow
\ker \alpha^q \rightarrow 0 \ $$ for any $q >0$. Combining with the long exact
sequence 
\begin{equation} \label{longp} \rightarrow \Hom(\sF, \sF)
  \stackrel{\alpha^0}{\longrightarrow} H^0(\sF) \rightarrow \Ext^1(\sI_{X/Y}
  \otimes \sL, \sF) \stackrel{p^1}{\longrightarrow} \Ext^1(\sF,
  \sF) \stackrel{\alpha^1}{\longrightarrow}
\end{equation} 
$$ H^1(\sF) \rightarrow
\Ext^2(\sI_{X/Y} \otimes \sL, \sF) \stackrel{p^2}{\longrightarrow} \Ext^2(\sF,
\sF) \stackrel{\alpha^2}{\longrightarrow} H^2(\sF)
$$
deduced from $0 \rightarrow \sO_Y \rightarrow \sF \rightarrow \sI_{X/Y}
\otimes \sL \rightarrow 0 \ $, we get the claim.

(i) From \eqref{longp} and the proven claim which leads to the fact that $p^1$
(resp. $p^2$) is the tangent map (resp. a map of obstruction spaces, mapping
obstructions to obstructions) of $p$, we get the smoothness of $p$ since $p^1$
is surjective and $p^2$ is injective. We will, however, give an independent
proof which one may use (slightly changed) to prove the remark below.
  
Let $(T,{\mathfrak m}_T) \rightarrow (S, {\mathfrak m}_S)$ be a small Artinian
surjection (i.e. of local Artinian $k$-algebras with residue fields $k$ whose
kernel ${\mathfrak a}$ satisfies ${\mathfrak a} \cdot {\mathfrak m}_T=0$). To
prove the (formal) smoothness of $p$, we must by definition show that the map
$${\rm Def}_{\sF,\sigma}(T) \rightarrow {\rm Def}_{\sF,\sigma}(S) \times_{{\rm
    Def}_{\sF}(S)} {\rm Def}_{\sF}(T)$$ is surjective. Let $\sigma_S:\sO_{Y
  \times S} \rightarrow \sF_S$ be a deformation of ${\sigma}$ to $S$ and let
$\sF_T$ be a deformation of $\sF_S$ to $T$. It suffices to find a map
${\sigma_T}: \sO_{Y \times T} \rightarrow \sF_T$ such that ${\sigma_T}
\otimes_T id_S = {\sigma_S}$, i.e. we must prove that $ H^0(\sF_T)
\rightarrow H^0(\sF_S)$ is surjective. Taking global sections of the short
exact sequence $$0 \rightarrow \sF \otimes_k {\mathfrak a} \simeq \sF_T
\otimes_T {\mathfrak a} \rightarrow \sF_T \rightarrow \sF_S \rightarrow 0 \ ,
$$ we get the surjectivity because $H^1(\sF) \otimes_k {\mathfrak a}=0$.

(ii) \ Again we have a long exact sequence
\begin{equation} \label{longq} \rightarrow \Ext^1(\sI_{X/Y} \otimes \sL,\sO_Y)
  \rightarrow \Ext^1(\sI_{X/Y} \otimes \sL, \sF)
  \stackrel{q^1}{\longrightarrow} \Ext^1(\sI_{X/Y} \otimes \sL, \sI_{X/Y}
  \otimes \sL) \rightarrow
\end{equation} 
$$ \Ext^2(\sI_{X/Y} \otimes \sL,\sO_Y) \rightarrow
\Ext^2(\sI_{X/Y} \otimes \sL,\sF) \stackrel{q^2}{\longrightarrow}
\Ext^2(\sI_{X/Y} \otimes \sL,\sI_{X/Y} \otimes \sL) \rightarrow
$$ containing maps $q^1$ (resp. $q^2$) which we may interpret as the tangent
map (resp. a map of obstruction spaces, which maps obstructions to
obstructions) of $q$. Indeed since ${\sE}xt^1(\sI_{X/Y},\sI_{X/Y}) \simeq
\sN_{X/Y}$ and ${\sH}om(\sI_{X/Y},\sI_{X/Y}) \simeq \sO_Y$, the assumption
$H^i(\sO_Y)=0$ for $i = 1, 2$ and the spectral sequence relating global and
local Ext-groups show $\Ext^1(\sI_{X/Y},\sI_{X/Y}) \simeq H^0(\sN_{X/Y})$ and
the injectivity of $\Ext^2(\sI_{X/Y},\sI_{X/Y}) \hookrightarrow
H^1(\sN_{X/Y})$ (see \cite{Vak}, the case $Y=\proj{3}$ was in fact proved in
\cite{K2}), as well as
\begin{equation} \label{skipl}
\Ext^i(\sI_{X/Y} \otimes \sL,\sI_{X/Y} \otimes \sL) \simeq
\Ext^i(\sI_{X/Y}, \sI_{X/Y} ) 
\end{equation}
for $i=1,2$. Hence we get the smoothness of
$q$ because $q^1$ is surjective and $q^2$ is injective by the assumption
$\Ext^2(\sI_{X/Y} \otimes \sL,\sO_Y) \simeq \Ext^2( \sF,\sO_Y)=0$. We will,
however, again give an independent proof using the definition of smoothness.

Let $T \rightarrow S$, ${\mathfrak a}$ and $\sigma_S:\sO_{Y \times S}
\rightarrow \sF_S$ be as in the proof of (i) above. Let $\sG_S = \coker
{\sigma_S}$ and let $\sG_T$ be a deformation of $\sG_S$ to $T$. By the theory
of extensions it suffices to show that the natural map
$$\Ext^1(\sG_{T},\sO_{Y \times T}) \rightarrow
\Ext^1(\sG_{S} ,\sO_{Y \times S}) $$ is surjective. Modulo isomorphisms we
refind this map in the middle of the long exact sequence
\begin{equation} \rightarrow \Ext^1(\sG_{T},\sO_{Y \times T} \otimes_T
  {\mathfrak a}) \rightarrow \Ext^1(\sG_{T},\sO_{Y \times T}) \rightarrow
  \Ext^1(\sG_{T}, \sO_{Y \times S}) \rightarrow \Ext^2(\sG_{T},\sO_{Y \times
    T} \otimes_T {\mathfrak a})
\end{equation}
Since $ \Ext^2(\sG_{T},\sO_{Y \times T} \otimes_T {\mathfrak a}) \simeq
\Ext^2(\sI_{X/Y} \otimes \sL,\sO_{Y}) \otimes_k {\mathfrak a} =0$ by
assumption we get the smoothness.

To see that $ {\rm Def}_{\sF,\sigma}(S) \simeq {\rm Def}_{X/Y,\xi}(S)$ are
isomorphic, take a deformation $\sO_{Y \times S}
\stackrel{\sigma_S}{\longrightarrow} \sF_S$ of $\sO_{Y}
\stackrel{\sigma}{\longrightarrow} \sF$. Since $\sF_S$ is flat, so are $\coker
\sigma_S$ and $(\coker \sigma_S) \otimes_{\sO_{Y_S}}(\sO_{Y_S} \otimes_{\sO_Y}
{\sL}^{-1})$. The former fits into a short exact sequence starting with $0
\rightarrow \sO_{Y \times S} \stackrel{\sigma_S}{\longrightarrow} \sF_S$, i.e.
we get an extension $\xi_S$ satisfying $\xi_S \otimes_S k = \xi$. The latter
is a flat deformation of $\sI_{X/Y}$. Thanks to the isomorphism
$\Ext^1(\sI_{X/Y},\sI_{X/Y}) \simeq H^0(\sN_{X/Y})$ and the injectivity
$\Ext^2(\sI_{X/Y},\sI_{X/Y}) \hookrightarrow H^1(\sN_{X/Y})$ above, one
knows that a deformation of $\sI_{X/Y}$ defines a deformation of $X$ in $ Y $,
i.e. we get an element $(X_S \subset Y_S) \in Hilb_{X/Y}(S)$ and hence we get
$ (X_S \subset Y_S, \xi_S) \in {\rm Def}_{X/Y,\xi}(S)$. This defines a map $
{\rm Def}_{\sF,\sigma}(S) \rightarrow {\rm Def}_{X/Y,\xi}(S)$. Since the
morphism the other way is just an obvious forgetful map, we get a functorial
isomorphism $ {\rm Def}_{\sF,\sigma}(S) \simeq {\rm Def}_{X/Y,\xi}(S)$, as
claimed in the Theorem.

{\rm (b)} To prove the first dimension formula, we 
continue \eqref{longp} to the left. Using $H^1(\sF)=0$ we get
$$\sum_{i=0}^1(-1)^{i+1} \ext^i(\sI_{X/Y} \otimes \sL, \sF)= \ext^1(\sF,
\sF) - hom(\sF, \sF)+h^0(\sF)$$ while \eqref{longq} (continued), $\Ext^2(
\sF,\sO_Y)=0$, ${\sH}om(\sI_{X/Y},\sI_{X/Y}) \simeq \sO_Y$ and \eqref{skipl}
show $$\sum_{i=0}^1(-1)^{i+1} \ext^i(\sI_{X/Y} \otimes \sL, \sF)=
\ext^1(\sI_{X/Y},\sI_{X/Y})-1 + \sum_{i=0}^1(-1)^{i+1} \ext^i(\sI_{X/Y}
\otimes \sL, \sO_Y).$$ Since ${\sE}xt^1(\sI_{X/Y},\sI_{X/Y}) \simeq \sN_{X/Y}$
it remains to show $ \sum_{i=0}^1(-1)^{i+1} \ext^i(\sI_{X/Y} \otimes \sL,
\sO_Y) = h^0(\omega_X \otimes \omega_Y^{-1}\otimes \sL^{-1})- \sum_{i=0}^2
(-1)^{i} h^i( {\sL}^{-1})$. Since ${\sH}om(\sI_{X/Y} \otimes \sL, \sO_Y) \simeq
\sL^{-1}$, ${\sE}xt^1(\sI_{X/Y}, \sO_Y) \simeq \omega_X \otimes \omega_Y^{-1}$
and $ \Ext^2(\sI_{X/Y} \otimes \sL, \sO_Y)=0$ (Remark~\ref{remver2}), we get
$\hom(\sI_{X/Y} \otimes \sL, \sO_Y)= h^0( {\sL}^{-1})$ and
$$\ext^1(\sI_{X/Y} \otimes \sL, \sO_Y)= h^0(\omega_X \otimes
\omega_Y^{-1}\otimes \sL^{-1}) + h^1( {\sL}^{-1})- h^2( {\sL}^{-1})$$ by the
spectral sequence relating global and local Ext-groups, and we get the first
dimension formula. 

Finally to see the last dimension formula and the genericness property, let $U
\subset {\rm Hilb}^{p}(Y)$ be a small enough open smooth connected subscheme
containing $(X)$ and let $\sI_{X_U/Y_U}$ be the sheaf ideal of $X_U \subset
Y_U:=Y \times U$, the universal object of $ {\rm Hilb}^{p}(Y)$ restricted to
$U$. Let $\sL_U:=\sL \otimes_{\sO_Y} \sO_{Y \times U} $. Looking to
\eqref{longq} which takes the form
\begin{equation} \label{fiblongq}0 \rightarrow H^0(\sO_Y) \rightarrow
  \Ext^1(\sI_{X/Y} \otimes \sL,\sO_Y) \rightarrow \Ext^1(\sI_{X/Y} \otimes
  \sL, \sF) \stackrel{q^1}{\longrightarrow} \Ext^1(\sI_{X/Y}, \sI_{X/Y})
  \rightarrow 0
\end{equation} 
and recalling that $q^1$ is the tangent map of $q: {\rm Def}_{\sF,\sigma}
\rightarrow {\rm Hilb}_{X/Y}$ and that $q$ is smooth, we can look upon the
fiber of $q$ as $ \Ext^1(\sI_{X/Y} \otimes \sL,\sO_Y)/k$. In the same way
since $\sF$ is stable and hence simple we can use the exact sequence
\eqref{longp} to see that the fiber of $p$ is isomorphic to $H^0(\sF)/k$.
Hence we get the second dimension formula since the functor $ {\rm Def}_{\sF}$
is prorepresented by the completion of the local ring of ${\rm M_Y}(P)$ at
$(\sF)$ (\cite{HL}, Thm. 4.5.1). More precisely the family $D:= \proj{}(
\Ext^1(\sI_{X_U/Y_U} \otimes \sL_U,\sO_{Y_U})^{\vee}) \rightarrow U$
parametrizes exactly extensions as in \eqref{serre} over $U$ and the
definition of a moduli space implies the existence of a morphism ${\overline
  p} : D \rightarrow {\rm M_{Y}}(P)$ whose corresponding local homomorphism at
$(\sF, \sigma)$ and $(\sF)$ induces $p$. Note that ${\overline p}$ is smooth
at$(X \subset Y,\xi)$ and hence maps the generic point of $D$ onto a generic
point of ${\rm M_Y}(P)$. This also proves that $\sF$ is a generic sheaf of
some component of ${\rm M_Y}(P)$ if and only if $X$ is generic in a
corresponding component of ${\rm Hilb}^{p}(Y)$ and we get all conclusions of
in the Theorem.
\end{proof}

\begin{remark} \label{graddefF} {\rm (a)} Suppose $Y$ is ACM and let
  $B:=H_{*}^0(\sO_Y)$. Applying $ H_{*}^0(-)$ onto \eqref{serre} we get an
  exact sequence $$0 \rightarrow B \stackrel{H_{*}^0(\sigma)}{\longrightarrow}
  F \rightarrow coker({H_{*}^0(\sigma))} \rightarrow 0 \ $$ inducing a long
  exact sequence (*) as in \eqref{longp} in which we have replaced the global
  ${\rm Ext}$-groups of sheaves with the corresponding graded ${\rm
    _0Ext}$-groups. Similar to ${\rm Def}_{\sF}$ (resp. ${\rm
    Def}_{\sF,\sigma}$) we may define local deformation functors ${\rm
    Def}_{F}$ (resp. ${\rm Def}_{F,H_{*}^0(\sigma)}$) on $ \underline{l}$ of
  flat graded deformations $F_S$ of $F$ (resp. $B \otimes_k S
  \stackrel{H_{*}^0(\sigma)}{\longrightarrow} F_S$ of $B
  \stackrel{H_{*}^0(\sigma)}{\longrightarrow} F$). There is a natural
  forgetful map $p_0: {\rm Def}_{F,H_{*}^0(\sigma)} \rightarrow {\rm Def}_{F}$
  whose tangent map fits into (*) and corresponds to $p^1$ in \eqref{longp}.
  Since ${\rm _0Ext}^1(B,F)=0$ in (*), it follows that $$p_0: {\rm
    Def}_{F,H_{*}^0(\sigma)} \rightarrow {\rm Def}_{F}$$ is smooth by the
  first (i.e. the cohomological) proof of Theorem~\ref{localserree} {\rm (i)}
  above.

  {\rm (b)} Suppose $Y$ is ACM and $ {_0\!\Hom_B}(F , M) = 0.$ Then we
  {\rm claim} that ${\rm Def}_{\sF} \simeq {\rm Def}_{F}$. Indeed by
  \eqref{seq} 
\begin{equation*}
  0 \rightarrow {_0\!\Ext_B^1 }(F , F) \rightarrow \Ext_{{\mathcal
      O}_{Y}}^1 (\sF ,\sF) \rightarrow 
  {_0\!\Ext_{\mathfrak m}^{2} }(F , F ) 
  \rightarrow  {_0\!\Ext_B^2 }(F , F) \rightarrow \Ext_{{\mathcal
      O}_{Y}}^2 (\sF ,\sF)
\end{equation*} 
is exact and ${_0\!\Ext_{\mathfrak m}^{2} }(F , F ) = {_0\!\Hom_B}(F , M)$ by
\eqref{mspect}. Hence we get the claim by the same cohomological argument as
used in Theorem~\ref{localserree} {\rm (i)} above. In the same way (or
directly) we can prove that ${\rm Def}_{\sF,\sigma} \simeq {\rm
  Def}_{F,H_{*}^0(\sigma)}$. It follows that the morphism $ p: {\rm
  Def}_{\sF,\sigma} \rightarrow {\rm Def}_{\sF}$ of Theorem~\ref{localserree}
is smooth.
\end{remark}

\section{Reflexive sheaves on $\proj{3}$ of small diameter}
As an application we concentrate on ${\rm M_{Y}}(c_1,c_2,c_3)$ with $Y =
\proj{3}$. Recalling the notion $F = H_{*}^0(\sF):= \oplus H^0(\sF(v))$, $M =
H_{*}^1(\sF)$ and $E = H_{*}^2(\sF)$, we have
\begin{theorem} \label{introrefl1} Let $\sF$ be a reflexive sheaf
  of rank $2$ on
  $\proj{3}$, 
  and suppose either 

  $\ \ \ \ (i) \ \ \ \ \ {_v\!\Hom_R}(F , M) = 0 \ \ \ {\rm for} \ v = 0 \ {\rm
    and} \ v = -4 \ , \ \ {\rm or}$ \\[-4mm]
  
  $\ \ \ \ (ii) \ \ \ \ {_v\!\Hom_R}(M , E) = 0 \ \ {\rm for} \ v = 0 \ {\rm
    and} \ v = -4 \ , \ \ {\rm or}$ \\[-4mm] 
  
  $\ \ \ \ (iii) \ \ \ {_0\!\Hom_R}(F , M) = 0 \ , \ \ {_0\!\Hom_R}(M , E) = 0
  $
  and that \\[-4mm]

   \ \ \ \ \ \ \ \ \ \ \ \ M is unobstructed as a graded module $( e.g. \
   {_{0}\!\Ext_R^2}(M ,M ) = 0 ) .$
   \\[2mm]
   Then $\sF$ is unobstructed. Moreover if ${_{0}\!\Ext_R^i}(M,M) = 0$ for $ i
   \geq 2$ and $\sF$ is stable, then ${\rm M}_{\proj{3}}(c_1,c_2,c_3)$ is
   smooth at $(\sF)$ and its dimension at $(\sF)$ is
\begin{equation*}
  \dim_{(\sF)}{\rm M}_{\proj{3}}(c_1,c_2,c_3) = {\rm ed}(\sF) +
  {_{-4}\!\hom_R}(F,M) +  {_{-4}\!\hom_R}(M , E)  + {_0\!\hom_R}(F,E) \ . 
\end{equation*}
Furthermore $ {_{-4}\!\hom_R}(M , E)= {_{0}\!\ext_R^1}(F,M)$ and
${_{-4}\!\hom_R}(M , E)+ {_0\!\hom_R}(F,E)= {_{-4}\!\ext_R^1}(F,F). $
\end{theorem}

Note that Theorem~\ref{introrefl1} applies to prove unobstructedness if $M=0$
(this case is known by \cite{Mir2}). The natural application of
Theorem~\ref{introrefl1} is to sheaves whose graded modules $M$ are
concentrated in a few degrees, e.g. $\diam M \le 2$. For such modules we can
prove more, namely that the sufficient conditions of unobstructedness of
Theorem~\ref{introrefl1} are quite close to being necessary conditions. Indeed
if the diameter of $M$ is one, they are necessary! Moreover in such cases a
minimal resolution of $\sF$ is often sufficient for computing the
$\Hom$-groups in the theorem (see also Lemma~\ref{smalldelta} below).

To find necessary conditions we consider the {\it cup product} or more
precisely its ``images'' in $ {_0\!\Hom_R}(F,E)$,
${_{-4}\!\Hom_R}(F,M)^{\vee}$ and ${_{-4}\!\Hom_R}(M,E)^{\vee}$ via some
natural maps, cf.\! \cite{W1}, \cite{F}, \cite{krao} and \cite{L2}, $ §
2$. 
Here we only include the cup product factorization given by (a) and hence
(b)(i) below, for which there is a proof in \cite{krao}, Prop.\! 3.6 of the
corresponding result for curves using Walter's factorization of $\alpha$ in
\eqref{mseq}. We remark that this result for curves, to our knowledge now, was
first proved by Fløystad (an easy consequence of Prop.\! 2.13 of \cite{F}).
For similarly generalizing the cases (ii) and (iii) of (b) we refer to
{\cite{krao}, Prop.\! 3.8. 
  Note that the necessary conditions in $(b)$ apply to many other sheaves than
  to those of diameter one (i.e. those with $M'=0$), e.g. they apply to
  Buchsbaum sheaves and to sheaves obtained by liaison addition of curves.

\begin {proposition} \label{introcupprod1} Let $\sF$ be a reflexive sheaf of
  rank $2$ on
  $\proj{3}$ and suppose ${_0\!\Ext_R^2}(M , M) = 0$. \\[2mm]
  (a) If the natural morphism
\begin{equation*}
  {_0\!\Hom_R}(F,M) \times {_0\!\Hom_R}(M , E) \longrightarrow  {_0\!\Hom_R}(F
  ,  E) 
\end{equation*}
(given by the composition) is non-zero, then $\sF$ is obstructed. \\[2mm]
(b) Suppose $M$ admits a decomposition $M=M' \oplus M_{[t]}$ as $R$-{\it
  modules} where the diameter of $ M_{[t]}$ is one and supported in degree
$t$. Then $\sF$ is obstructed
provided  \\[-3mm]

$  \ \  (i) \ \ \ 
{_0\!\Hom_R}(F, M_{[t]}) \neq 0 \ \ {\rm and} \ \ \ {_0\!\Hom_R}( M_{[t]}, E)
\neq 0 \ , \ \ or $  \\[-3mm]

$  \ \ (ii) \ \ 
{_{-4}\!\Hom_R}(F,M_{[t]}) \neq 0 \ \ {\rm and} \ \ \ {_0\!\Hom_R}(M_{[t]},
E) \neq 0  \ , \ \ or $  \\[-3mm]

$ \ \ (iii) \ \ 
{_0\!\Hom_R}(F,M_{[t]}) \neq 0 \ \ {\rm and} \ \ \ {_{-4}\!\Hom_R}(M_{[t]}, E)
\neq 0 \ . $
\end{proposition}

\begin{proof}[Proof (of Theorem~\ref{introrefl1} and
  Proposition~\ref{introcupprod1})] 

  To prove the results we may replace $\sF$ by $\sF(j)$ for $j>>0$ because,
  for both results, the assumptions as well as the conclusions 
  hold for $\sF$ iff they hold for $\sF(j)$. In particular we may assume $
  H^1(\sF(v))=0$ for $v \le 0$ and hence $\Ext^2(\sF,\sO_{\proj{3}})^{\vee}
  \simeq H^1(\sF(-4))=0$. It follows that the maps $p$ and $q$ of
  Theorem~\ref{localserree} are smooth. From the Hartshorne-Serre
  correspondence we get an exact sequence $$0 \rightarrow R \rightarrow F
  \rightarrow I_{X}(c_1) \rightarrow 0 \ $$ which implies $M = H_{*}^1(\sF)
  \simeq H_{*}^1(\sI_X(c_1))$. We also get the exact sequence $$0 \rightarrow
  E = H_{*}^2(\sF) \rightarrow H_{*}^2(\sI_X(c_1)) \rightarrow
  H_{*}^3(\sO_{\proj{3}}) . $$ Using these sequences and $
  H^1(\sI_X(c_1+v))=0$ for $v \le 0$ we get
  \begin{equation} \label{comp}{_v\!\Hom_R}(F,M) \simeq {_v\!\Hom_R}(I_X,
    H_{*}^1(\sI_X))\ \ {\rm \ and \ } \ \ {_v\!\Hom_R}(M,E) \simeq
    {_v\!\Hom_R}( 
    H_{*}^1(\sI_X), H_{*}^2(\sI_X))
  \end{equation} for $-4 \le v \le 0$ because ${_v\!\Hom_R}(R,M)=0$ and $
  {_v\!\Hom_R}(M, H_{*}^3(\sO_{\proj{3}})) \simeq \ _{v}\!\Ext_{\mathfrak
    m}^{4}(M,R) \simeq M_{-v-4}^{\vee}=0$ by \eqref{mspect} and
  \eqref{mduality}. Now recall that  
  we in \cite{krao} proved results similar to Theorem~\ref{introrefl1} and
  Proposition~\ref{introcupprod1} for the unobstructedness (resp.
  obstructedness) of $X$ with the
  difference that the $\Hom$-groups, $ {_v\!\Hom_R}(H_{*}^i(\sF),
  H_{*}^{i+1}(\sF))$ for $\sF$ 
  were exchanged by the corresponding groups,  $
  {_v\!\Hom_R}(H_{*}^i(\sI_X), H_{*}^{i+1}(\sI_X))$  for $\sI_X$.
  Therefore  \eqref{comp} and Theorem~\ref{localserree} show that $\sF$ is
  unobstructed in Theorem~\ref{introrefl1} (resp. obstructed in
  Proposition~\ref{introcupprod1}) because $X$ is correspondingly unobstructed
  (resp. obstructed) by the results of \cite{krao},  Thm.\! 2.6 (resp. Prop.\!
  3.6 and Thm.\! 3.2) and  Remark~\ref{extend} below. 

  To prove the dimension formula we suppose ${_0\!\Ext_R^i}(M,M)=0$ for $2 \le
  i \le 4$. With this assumption the map $\alpha$ in \eqref{mseq} is zero by
  Walter's observation. Note that there is a corresponding connecting map
  $\alpha(N_1,N_2):\ _{0}\!\Ext_{\mathfrak m}^{2}(N_1,N_2) \rightarrow \
  _0\!\Ext_R^{2}(N_1,N_2)$ appearing in \eqref{seq}. Indeed
  $\alpha=\alpha(I_X,I_X)$. We {\it claim} that $\alpha(F,F)=0$. To see it we
  use the functoriality of the sequence \eqref{seq} and $\alpha=0$. Since the
  natural map $ {_0\!\Ext_R^2}(I_X(c_1),F) \rightarrow
  {_0\!\Ext_R^2}(I_X(c_1),I_X(c_1))$ is an isomorphism by $$
  {_0\!\Ext_R^2}(I_X(c_1),F)^{\vee} \simeq \ _{-4}\!\Ext_{\mathfrak
    m}^{2}(F,I_X(c_1)) \simeq {_{-4}\!\Hom_R}(F,M),$$ $$
  {_0\!\Ext_R^2}(I_X(c_1),I_X(c_1))^{\vee} \simeq \ _{-4}\!\Ext_{\mathfrak
    m}^{2}(I_X,I_X) \simeq {_{-4}\!\Hom_R}(I_X, H_{*}^1(\sI_X))$$ and
  \eqref{comp}, we get $\alpha(I_X(c_1),F)=0$. In the same way the natural map
  $ {_0\!\Ext_{\mathfrak m}^2}(I_X(c_1),F) \rightarrow \ _{0}\!\Ext_{\mathfrak
    m}^{2}(F,F)$ is an isomorphism by \eqref{seq} (i.e. both groups are
  naturally dual to $ {_{-4}\!\Hom_R}(F,M)$ by \eqref{mduality}) and we get
  the claim from $\alpha(I_X(c_1),F)=0$. Now using the fact that the
  projective dimension of $F$ is $2$, the proven claim and \eqref{seq}, we get
  an exact sequence
\begin{equation} \label{exFF}
  0 \rightarrow {_0\!\Ext_R^2}(F,F) \rightarrow  \Ext_{{\mathcal
      O}_{\proj{}}}^2 (\sF,\sF) \rightarrow  {_0\!\Ext_{\mathfrak m}^{3}}(F,F)
  \rightarrow  0 .
\end{equation}
As above ${_0\!\Ext_R^2}(F,F)^{\vee} \simeq {_{-4}\!\Hom_R}(F,M)$ and
similarly ${_0\!\Ext_{\mathfrak m}^{3}}(F,F)^{\vee} \simeq
{_{-4}\!\Ext_R^1}(F,F)$ by \eqref{mduality}. We get $ \ext_{{\mathcal
    O}_{\proj{}}}^2 (\sF,\sF) = {_{-4}\!\hom_R}(F,M)+{_{-4}\!\ext_R^1}(F,F).$
Using \eqref{mspect} we get an exact sequence
\begin{equation} \label{exmFF} 0 \rightarrow {_0\!\Ext_R^1}(F,M) \rightarrow
  {_0\!\Ext_{\mathfrak m}^{3}}(F,F) \rightarrow {_0\!\Hom_R}(F,E) \rightarrow
  {_0\!\Ext_R^2}(F,M) \rightarrow
\end{equation}
and hence ${_{-4}\!\ext_R^1}(F,F) = {_{0}\!\ext_R^1}(F,M) + {_0\!\hom_R}(F,E)
$ because
\begin{equation} \label{m1}
{_{0}\!\Ext_R^2}(F,M) \simeq \ _{-4}\!\Ext_{\mathfrak
  m}^{2}(M,F)^{\vee} \simeq {_{-4}\!\Hom_R}(M,M)^{\vee} \simeq \
_{0}\!\Ext_{\mathfrak m}^{4}(M,M) \simeq {_0\!\Ext_R^4}(M,M)=0 .
\end{equation} 
By the arguments of \eqref{m1} we also get ${_{-4}\!\Ext_R^1}(M,M)^{\vee}
\simeq \ {_0\!\Ext_R^3}(M,M)=0$ and $${_{0}\!\ext_R^1}(F,M) = \
_{-4}\!\ext_{\mathfrak m}^{3}(M,F) = {_{-4}\!\hom_R}(M , E)$$ and putting
things together we are done.
\end{proof}

\begin{remark} \label{extend} Thm.\! 2.6(iii) of \cite{krao} actually proves a
  slightly weaker statement than needed to prove
  Theorem~\ref{introrefl1}(iii). However, putting different results of e.g.
  \cite{krao} together we get what we want. Indeed we claim that a curve $X
  \subset \proj{3}$ is unobstructed provided
  \begin{equation} \label{cc} {_0\!\Hom_R}(I_X,H_{*}^1(\sI_X))= 0 \ , \ \
    {_0\!\Hom_R}(H_{*}^1(\sI_X), H_{*}^{2}(\sI_X)) = 0 \ ,
  \end{equation}   {\rm and } $H_{*}^1(\sI_X)$ {\rm is \ unobstructed \ as \
    a \ graded  \ module \ 
    (e.g.} \ ${_{0}\!\Ext_R^2}(H_{*}^1(\sI_X),H_{*}^1(\sI_X)) = 0${\rm )}.
  This is mainly a consequence of results proven in \cite{MDP1} by
  Martin-Deschamps and Perrin. Indeed their smoothness theorem for the
  morphism from the Hilbert scheme of constant cohomology, $H(d,g)_{cc}$, onto
  the scheme of ``Rao modules'' (Thm.1.5, p.\! 135) combined with their
  tangent space descriptions (pp. 155-156), or more precisely combined with
  Prop.\! 2.10 of \cite{krao} which states that the vanishing of the two
  $\Hom$-groups in \eqref{cc} leads to an isomorphism $H(d,g)_{cc} \simeq
  H(d,g)$ at $(X)$, we conclude easily.
\end{remark}

We can compute the number  ${_0\!\hom_R}(F,E)$ in terms of the
graded Betti numbers  $\beta_{j,i}$ of $F$;
\begin{equation} \label{resoluMF}
0 \rightarrow \bigoplus_i
R(-i)^{\beta_{3,i}} \rightarrow \bigoplus_i R(-i)^{\beta_{2,i}} \rightarrow
\bigoplus_i R(-i)^{\beta_{1,i}} \rightarrow F \rightarrow 0 \
\end{equation} (sheafifying, we get the ``resolution'' of $\sF$ in the
introduction),  by using the following result.

\begin {lemma} \label{smalldelta} Let $\sF$ be a reflexive sheaf of
  rank $2$ on $\proj{3}$, and suppose $ {_{-4}\!\Hom_R}(F,F)=0$. Then

$$ {_0\!\hom_R}(F,E) = \sum _{i}(\beta_{1,i}
- \beta_{2,i} + \beta_{3,i}) \cdot (h^2(\sF(i))-h^3(\sF(i))) \ . \ $$
\end{lemma}

\begin{proof} Recall $E:= H_{*}^2(\sF) \simeq H_{\mathfrak m}^{3}(F)$. If we
  apply ${_v\!\Hom_R}(- ,E)$ to the minimal resolution \eqref{resoluMF} we
  get a complex
  \begin{equation} \label{resolutMF} 0 \rightarrow {_0\!\Hom_R}(F,
    H_{*}^2(\sF)) \rightarrow \bigoplus_i H^2(\sF(i))^{\beta_{1,i}}
    \rightarrow \bigoplus_i H^2(\sF(i))^{\beta_{2,i}} \rightarrow
    \bigoplus_i H^2(\sF(i))^{\beta_{3,i}} \rightarrow 0 \ .
  \end{equation} 
  Since the alternating sum of the dimension of the terms in a complex equals
  the alternating sum of the dimension of its homology groups, it suffices to
  show ${_0\!\Ext}_R^1(F,E)=0$, ${_0\!\Ext}_R^2(F,E) \simeq
  {_0\!\Hom_R}(F,H_{*}^3(\sF))$ and that 
  \begin{equation} \label{HomH3} {_0\!\hom_R}(F,H_{*}^3(\sF)) = \sum
    _{i}(\beta_{1,i} - \beta_{2,i} + \beta_{3,i}) \cdot h^3(\sF(i)).
  \end{equation}
  Using \eqref{mspect} and that $ _{0}\!\Ext_{\mathfrak m}^{4}(F,F) \simeq
  {_{-4}\!\Hom_R}(F,F)^{\vee} =0$ by assumption, we get $
  {_0\!\Ext_R^1}(F,H_{\mathfrak m}^{3}(F))=0$ and an exact sequence $$0
  \rightarrow {_{0}\!\Hom_R}(F,H_{\mathfrak m}^{4}(F)) \rightarrow
  {_0\!\Ext_R^2}(F,H_{\mathfrak m}^{3}(F)) \rightarrow \ _0\!\Ext_{\mathfrak
    m}^{5}(F,F) \rightarrow {_{0}\!\Ext_R^1}(F,H_{\mathfrak m}^{4}(F))
  \rightarrow 0 . $$ Since we have $_{0}\!\Ext_{\mathfrak m}^{5}(F,F)=0$ the
  proof is complete provided we can show \eqref{HomH3}. To show it, it is
  sufficient to see that \eqref{resolutMF}, with $H^2$ replaced by $H^3$,
  is exact. Since we have $ _0\!\Ext_{\mathfrak m}^{i}(F,F)=0$ for $i=5,6$ by
  duality, we get $ {_{0}\!\Ext_R^i}(F,H_{\mathfrak m}^{4}(F))= 0 $ for
  $i=1,2$ by \eqref{mspect} and we are done.
\end{proof}

\begin{remark} \label{lowbound} For later use we remark that if we apply
  ${_0\!\hom_R}(-,M)$, $M = H_{*}^1(\sF)$ to \eqref{resoluMF} we get
  \begin{equation} \label{lemm}  \sum_{ i = 0 }^{2} (-1) ^{i}
    ~ {_0\!\ext ^i} (F , M) =  \sum _{i}(\beta_{1,i}
- \beta_{2,i} + \beta_{3,i}) \cdot h^1(\sF(i)) \ ,
\end{equation} 
cf. \eqref{resolutMF}. Suppose $\sF$ is reflexive and $ {_{-4}\!\Hom_R}(F,F)
=0$. Using \eqref{mspect} as in \eqref{exmFF} and the proof above we get
$\sum_{ i = 2 }^{3} (-1)^{i} \, {_0\!\ext_{\mathfrak m}^i} (F , F)=\sum_{ i =
  0 }^{2} (-1) ^{i} ~ {_0\!\ext ^i} (F , M) - {_0\!\hom_R}(F,E)$. Hence we
have $$\sum_{ i = 2 }^{3} (-1)^{i} \, {_0\!\ext_{\mathfrak m}^i} (F , F)= \sum
_{i}(\beta_{1,i} - \beta_{2,i} + \beta_{3,i}) \cdot (
h^1(\sF(i))-h^2(\sF(i))+h^3(\sF(i))).$$
\end{remark}
It is easy to substitute the non-vanishing of the $\Hom$-groups of
Theorem~\ref{introrefl1} by the non-triviality of certain graded Betti numbers
in the minimal resolution of $F$. Indeed we have

\begin{theorem} \label{refl3} Let $\sF$ be a reflexive sheaf of rank $2$
  on $\proj{3}$ and suppose $M \neq 0$ is of diameter
  $1$ and concentrated in degree $c$. 
  Then $\sF$ is obstructed if and only if
  $$
  \beta_{1,c} \cdot \beta_{2,c+4} \neq 0 \ \ \ {\rm or} \ \ \ \beta_{1,c+4}
  \cdot \beta_{2,c+4} \neq 0 \ \ \ {\rm or} \ \ \ \beta_{1,c} \cdot
  \beta_{2,c}\neq 0 \ . $$ Moreover if $\sF$ is an unobstructed stable sheaf
  and $M$ is $r$-dimensional, then the dimension of the moduli scheme
  ${\rm M}_{\proj{3}}(c_1,c_2,c_3)$ at $(\sF)$ is
\begin{equation*}
  \dim_{(\sF)}{\rm M}_{\proj{3}}(c_1,c_2,c_3) =  {\rm ed}(\sF)+
  {_0\!\hom_R}(F,E) + 
  r(\beta_{1,c+4} +  \beta_{2,c}). 
\end{equation*}
\end{theorem}

Before proving Theorem~\ref{refl3}, we remark that we have the following
result

\begin {proposition} \label{propreflobstr} Let $\sF$ be a reflexive sheaf of
  rank $2$ on $\proj{3}$ and suppose $M \neq 0$ is of diameter $1$. Then $\sF$
  is obstructed iff (at least) one of
  the following conditions hold  \\[-3mm]

$\ \ \ \ (a) \ \ \ \ \ {_0\!\Hom_R}(F , M) \neq 0 \ \ \ {\rm and} \ \ \
{_0\!\Hom_R}(M , E) \neq 0 \ ,  $ \\[-3mm]

 $\ \ \ \ (b) \ \ \ \ {_{-4}\!\Hom_R}(F , M) \neq 0 \ \ \ {\rm and} \ \ \
  {_0\!\Hom_R}(M , E) \neq 0  \ ,  $  \\[-3mm]

  $\ \ \ \ (c) \ \ \ \ \ {_0\!\Hom_R}(F , M) \neq 0 \ \ \ \ {\rm and} \ \
  {_{-4}\!\Hom_R}(M , E) \neq 0 \ . $
\end {proposition}

Indeed if $\sF$ is unobstructed, then it is a simple reformulation of
Theorem~\ref{introrefl1} to see that we have either (a) or (b) or (c). The
converse follows immediately from  Proposition~\ref{introcupprod1} by letting  
$M'=0$.

\begin{proof}[Proof (of Theorem~\ref{refl3})] 
  By applying ${_v\!\Hom_R}(- ,M)$ to the minimal resolution \eqref{resoluMF}
  we get
\begin{equation} \label{corbetti1} \ {_0\!\hom_R}(F,M) = r \beta_{1,c} \ \
{\rm and} \ \ {_{-4}\!\hom_R}(F, M) = r \beta_{1,c+4} \ .
\end{equation}
because $\mathfrak m \cdot M = 0$ and $ \oplus_i R(-i)^{\beta_{3,i}} =
R(-c-4)^{r}$. Moreover ${_{-v-4}\!\Ext_R^1}(F,M)^{\vee} \simeq {_v\!\Hom_R}(M,
E)$. Interpreting ${_{-v-4}\!\Ext_R^1}(F ,M)$ similarly via the minimal
resolution \eqref{resoluMF} of $F$, we get
\begin{equation} \label{corbetti2} {_0\!\hom_R}(M, E) = r\beta_{2,c+4} \quad
  {\rm and} \quad {_{-4}\!\hom_R}(M, E) = r \beta_{2,c}\ .
\end{equation}
Since $r \neq 0$, we get the unobstructedness criterion and the dimension
formula of Theorem~\ref{refl3} from Proposition~\ref{propreflobstr} and
Theorem~\ref{introrefl1}.
\end{proof}

\begin{theorem} \label{diam1} Every irreducible component of ${\rm
    M}_{\proj{3}}(c_1,c_2,c_3)$ whose generic sheaf $\sF$ satisfies $\diam M
  \le 1$ is {\it reduced} (i.e. generically smooth).
\end{theorem}

\begin{proof} By replacing $\sF$ by $\sF(j)$ for $j>>0$ (cf. first part of the
  proof of Theorem~\ref{introrefl1}), we can use the Hartshorne-Serre
  correspondence to get a corresponding curve $X$ such that all assumptions of
  Theorem~\ref{localserre1} are satisfied. Hence $X$ is generic and $\diam
  H_{*}^1(\sI_X) \le 1$. Since it is proved in \cite{krao}, Cor.\! 4.3 that a
  generic curve $X$ of diameter at most one is unobstructed, it follows by
  Theorem~\ref{localserre1} that
  $\sF$ is unobstructed, i.e. the corresponding component
  ${\rm M}_{\proj{3}}(c_1,c_2,c_3)$ is generically smooth and we are done.
\end{proof}

 \begin{example} \label{ex1} Using some results of Chang on
   $\Omega$-resolutions of Buchsbaum curves (\cite{C} or \cite{W1}, Thm. 4.1),
   one shows that there exists a smooth connected curve $X$ of diameter 1
   satisfying $h^0(\sI_X(e)) = 1$, $h^1(\sI_X(e)) = r$, $h^1(\sO_X(e)) = b$,
   $h^1(\sO_X(v)) = 0$ for $v>e$
   and with $e = 1 + b + 2r$ and $\Omega$-resolution
\begin{equation}  \label{omeg}
 0 \rightarrow \sO_{\proj{}}(-2)^{3r-1} \oplus
    \sO_{\proj{}}(-4)^{b} \rightarrow \sO_{\proj{}} 
    \oplus \Omega^{r}
    \oplus \sO_{\proj{}}(-3)^{b-1} \rightarrow \sI_X(e)
    \rightarrow 0
\end{equation}
for any pair $(r,b)$ of positive integers (cf. \cite{krao}, Ex. 3.12).
Moreover the degree and genus of $X$ is $d = { e+4 \choose 2} - 3r - 7$ and $g
= (e+1)d - { e+4 \choose 3}+ 5$. Recalling that $\Omega$ corresponds to the
first syzygy in the Koszul resolution of the regular sequence $\{
X_0,X_1,X_2,X_3\} $, we get an exact sequence
\begin{equation*} \label{ex11}
 \ \ 0 \rightarrow
  \sO_{\proj{}}(-4) 
  \rightarrow \sO_{\proj{}}(-3)^{4} \rightarrow
  \sO_{\proj{}}(-2)^{6} \rightarrow  \Omega \rightarrow 0 
\end{equation*} 
Hence we can use the mapping cone construction to show that there is a
resolution
\begin{equation} \label{ex22} 0 \rightarrow \sO_{\proj{}}(-4)^{r} \rightarrow
  \sO_{\proj{}}(-4)^{b} \oplus \sO_{\proj{}}(-3)^{4r} \oplus
  \sO_{\proj{}}(-2)^{3r-1} \rightarrow \sO_{\proj{}}(-2)^{6r} \oplus
  \sO_{\proj{}} \oplus \sO_{\proj{}}(-3)^{b-1} \rightarrow \sI_{X}(e)
  \rightarrow 0 
\end{equation}
where we possibly may skip the factor $ \sO_{\proj{}}(-2)^{3r-1}$ (and reduce
$\sO_{\proj{}}(-2)^{6r}$ to $ \sO_{\proj{}}(-2)^{3r+1}$) to get a minimal
resolution. Instead of looking into this problem, we will illustrate
\cite{krao}, Thm.\! 4.1, which makes a deformation theoretic improvement to a
theorem of Rao (\cite{R} Thm.\! 2.5). Indeed since the composition of the
leftmost non-trivial map in \eqref{ex22} with the projection onto
$\sO_{\proj{}}(-2)^{3r-1}$ is zero (by Rao's theorem), there is by
\cite{krao}, Thm.\! 4.1 a deformation with constant cohomology and Rao module
to a curve $X'$ which makes $\sO_{\proj{}}(-2)^{3r-1}$ redundant (no matter
whether the original factor was redundant or not)! So we certainly may skip
the factor $ \sO_{\proj{}}(-2)^{3r-1}$ (and reduce $\sO_{\proj{}}(-2)^{6r}$ to
$ \sO_{\proj{}}(-2)^{3r+1}$), at least after a deformation.

Now, by the Hartshorne-Serre correspondence there is a reflexive sheaf $\sF$
given by
\begin{equation} \ 0 \rightarrow \sO_{\proj{}}
  \stackrel{\sigma}{\longrightarrow} \sF \rightarrow \sI_{X'}(e+4) \rightarrow
  0.
\end{equation} which combined with the Horseshoe lemma \cite{Wei} leads
to the following minimal resolution of $\sF$,
\begin{equation} \label{ex10} 0 \rightarrow \sO_{\proj{}}^{r} \rightarrow
  \sO_{\proj{}}^{b} \oplus \sO_{\proj{}}(1)^{4r} \rightarrow
  \sO_{\proj{}}(2)^{3r+1} \oplus \sO_{\proj{}}(4) \oplus
  \sO_{\proj{}}(1)^{b-1} \oplus \sO_{\proj{}}\rightarrow \sF
  \rightarrow 0 \ .
\end{equation}
Note that $h^1(\sF(-4)) = h^1(\sI_X(e)) = r$, i.e. the number $c$ of
Theorem~\ref{refl3} is $c=-4$. From \eqref{ex10} we see that
$\beta_{2,0}=b \ne 0$ and $\beta_{1,-4}=1$. By Theorem~\ref{refl3}, $\sF$
is obstructed.

Computing Chern classes $c_i$ of $\sF$ we get $c_1 = e+4$, $c_2=d= { c_1
  \choose 2} - 3r - 7$ and $c_3 = { c_1 \choose 3} - { c_1 \choose 2}(3r +
7)+6r+22$. The simplest case is $(r,b) = (1,1)$, which yields a reflexive
sheaf $\sF$ whose normalized sheaf $\sF(-4)$ is semistable and with Chern
classes $(c_1',c_2',c_3')$ = $(0,2,4)$ (the corresponding curve $X$ has $d=
18$, $g = 39$ and is Sernesi's example of an obstructed curve, \cite{Se} or
\cite{EF}, see also \cite{MT} which thoroughly studies ${\overline
  M}_{\proj{3}}(0,2,4)$ and \cite{Mir5} which uses Sernesi's example to show
the existence of a stable rank 3 obstructed {\rm vector bundle}). For $(r,b)
\ne (1,1)$, then $e > 4$ and we see easily that the sheaves constructed above
are stable. If $(r,b)=(2,1)$, then the normalized sheaf has Chern classes
$(c_1',c_2',c_3')$ = $(0,7,24)$ while $(r,b)=(1,2)$ yields stable sheaves with
$(c_1',c_2',c_3')$ = $(-1,6,22)$. One may show that all curves corresponding
to the sheaves of the case $(r,1)$ satisfy $h^1(\sN_X) = 1$. The ideal of the
local ring of $\HH(d,g)$ at $(X)$ is generated by a single element, which is
irreducible for $r>1$. Indeed the case $(r = 1)$ corresponds to a sheaf which
sits in the intersection of two irreducible components of ${\rm
  M}_{\proj{3}}(c_1,c_2,c_3)$, while for $r > 1$, the irreducibility of can be
used to see that $(\sF)$ belongs to a unique irreducible component of ${\rm
  M}_{\proj{3}}(c_1,c_2,c_3)$.
\end{example} 

In the examples~\ref{can} and \ref{mum} the diameter of $M$ of the obstructed
generic sheaves is $3$. Combining the results of this paper with the large
number of non-reduced components one may find in \cite{K1} we can easily
produce similar examples for every $\diam M \ge 3$. Indeed, as is well known,
a smooth cubic surface $X \subset \proj{3}$ satisfies $Pic(X) \simeq {\mathbb
  Z}^{\oplus 7}$. It follows from the main theorem of \cite{K1} (or of
\cite{N}) that the general curve which corresponds to $(3\alpha,\alpha^5,2)
\in {\mathbb Z}^{\oplus 7}$ is the generic curve of a non-reduced component of
$ \HH(d,g)$ for every $ \alpha \ge 4$ (Mumford's example corresponds to $
\alpha = 4$). The diameter is $2 \alpha -5$. In the same way the general curve
which corresponds to $(3\alpha+1,\alpha^5,2) \in {\mathbb Z}^{\oplus 7}$ is
the generic curve of a non-reduced component of $\HH(d,g)$ with $ \diam M = 2
\alpha -4$ for every $ \alpha \ge 4$. Using Theorem~\ref{localserre1} for
$c_1=2$ we get non-reduced components of ${\rm M}_{\proj{3}}(c_1,c_2,c_3)$ for
every $\diam M(\sF) \ge 3$, $\sF$ the generic sheaf. Thanks to
Theorem~\ref{diam1} there is only one case left and we expect:

\begin{conjecture} Every irreducible component of ${\rm
    M}_{\proj{3}}(c_1,c_2,c_3)$ whose generic sheaf $\sF$ satisfies $\diam M =
  2$ is {\it reduced} (i.e. generically smooth).
\end{conjecture}

There are some evidence to the conjecture, namely that {\it every} Buchsbaum
curve of diameter at most 2 admits a generization in $\HH(d,g)$ which is
unobstructed (\cite{krao}, Cor.\! 4.4), i.e. belongs to a generically smooth
irreducible component. By the arguments in the proof of Theorem~\ref{diam1}
{\it every} Buchsbaum sheaf of diameter at most 2 must belong to a generically
smooth irreducible component of some ${\rm M}_{\proj{3}}(c_1,c_2,c_3)$.

\section {A lower bound of $\dim {\rm M}_{\proj{3}}(c_1,c_2,c_3)$.}
In this section we want to give a lower bound of the dimension of any
irreducible component of ${\rm M}_{\proj{3}}(c_1,c_2,c_3)$ in terms of the
graded Betti numbers of a minimal resolution of the graded $R$-module $F =
H_{*}^0(\sF)$, see \eqref{resoluMF}. The lower bound is straightforward to
compute provided we know the dimension of the cohomology groups $ H^i(\sF(v))$
for any $i$ and $v$. It is well known that ${\rm ed}(\sF)=8c_2-2c_1^2-3$ is a
lower bound, but there are many examples of so-called oversized irreducible
components whose dimension is strictly greater that ${\rm ed}(\sF)$. Our lower
bound is usually much closer to the actual dimension of the oversized
components provided $H_{*}^1(\sF)$ is ``small''. If a component of ${\rm
  M}_{\proj{3}}(c_1,c_2,c_3)$ is generically smooth, we also include a formula
for dimension of the component which is a sum of the lower bound and a
correction number which we make explicit.


\begin{definition} \label{deltaF} 
  If $\sF$ is a reflexive sheaf on  $\proj{3}$, we let
\begin{equation*}
  \delta^{j} = \sum _{i} (\beta_{1,i} -\beta_{2,i}+ \beta_{3,i}) \cdot
  h^j(\sF(i)) .  
\end{equation*}
\end{definition}

\begin{remark} \label{explici} If $\sF$ is reflexive on $\proj{3}$ and ${
    _{-4}\!\Hom_R}(F,F)= 0$, then
  $${ _0\!\hom_R}(F, E) = \delta^{2}- \delta^{3}$$
  by Lemma~\ref{smalldelta}. This makes the dimension formulas of Theorems of
  \ref{introrefl1} and \ref{diam1} more explicit.
\end{remark}

\begin{proposition}\label{propeulerIF}
  Let $\sF$ be a reflexive sheaf on $\proj{3}$ satisfying ${
    _{-4}\!\Hom_R}(F,F)= 0$. Then
\begin{equation*}
  {_0\!\ext_R^1 }(F,F )- {_0\!\ext_R^2 }(F,F ) =  { _{0}\!\hom_R}(F,F) -
  \delta^0 = {\rm  ed}(\sF) +  \delta^2 - \delta^1 - \delta^3.
\end{equation*}
\end{proposition}

\begin {proof} To see the equality to the left, we apply ${_0\!\Hom_R}(-, F)$
  to the resolution \eqref{resoluMF}. We get
  \begin{equation*} { _0\!\hom_R}(F, F) -\ _0\!\ext_R^1(F , F) +
    {_0\!\ext_R^2}(F , F) = \delta^0 \ . \ 
\end{equation*}
Moreover the right hand equality follows from \eqref{mduality} and \eqref{seq}.
Indeed we have already considered the consequences of \eqref{mduality} in
Lemma~\ref{smalldelta} and Remark~\ref{lowbound}. We have
\begin{equation*} {_0\!\ext_m^2}(F, F) - {_0\!\ext_ m^3}(F , F)
  = \delta^1 - \delta^2 + \delta^3 \ 
\end{equation*}
by Remark~\ref{lowbound}. Combining with the exact sequence \eqref{mseq} which
implies
\begin{equation*} {\rm ed} (\sF) = {_0\!\ext_R^1 }(F,F )-
  {_0\!\ext_R^2 }(F,F ) + {_0\!\ext_m^2}(F , F) - {_0\!\ext_ m^3}(F , F) \ ,
\end{equation*}
we get the last equality.
\end{proof}

\begin{theorem}\label{thmeulerIF}
  Let $\sF$ be a stable reflexive sheaf on $\proj{3}$. Then the dimension of
  ${\rm M}_{\proj{3}}(c_1,c_2,c_3)$ at $(\sF)$ satisfies
\begin{equation*}
  \dim_{(\sF)}{\rm M}_{\proj{3}}(c_1,c_2,c_3) \ge  1 -
  \delta^0 = {\rm  ed}(\sF) +  \delta^2 - \delta^1 - \delta^3.
\end{equation*}
Moreover if $\sF$ is a generic sheaf of a generically smooth component $V$ of
${\rm M}_{\proj{3}}(c_1,c_2,c_3)$ and $M= H_{*}^1(\sF)$, then
$$\dim V = {\rm ed}(\sF) + \delta^2 - \delta^1
- \delta^3 +\ {_{-4}\!\hom_R}(F , M) $$ where ${_{-4}\!\Hom_R}(F , M)$ is the
kernel of the map $$ \bigoplus_i H^1(\sF(i-4))^{\beta_{1,i}} \longrightarrow
\bigoplus_i H^1(\sF(i-4))^{\beta_{2,i}} $$ induced by the corresponding map in
\eqref{resoluMF}.
\end{theorem}

\begin{remark} Let $\sF$ be a stable reflexive sheaf on $\proj{3}$. \\
  {\rm (i)} If $M=0$, then $\delta^1=0$ and we can use
  Theorem~\ref{introrefl1} and Remark~\ref{explici} to see that the lower
  bound of Theorem~\ref{thmeulerIF} is equal to $ \dim_{(\sF)}{\rm
    M}_{\proj{3}}(c_1,c_2,c_3)$. 
  This coincides with \cite{Mir2}. \\ 
  {\rm (ii)} If $\diam M=1$ and $\sF$ is a generic unobstructed sheaf, then
  the lower bound coincides with ${\rm ed}(\sF)+ {_0\!\hom_R}(F,E) + r
  \beta_{2,c}$ of Theorem~\ref{refl3} because $ r \beta_{1,c}=0$ for a generic
  sheaf by \cite{krao}, Cor.\! 4.4 and the proof of Theorem~\ref{diam1}.
  Moreover in this case the correction number ${_{-4}\!\hom_R}(F , M)$ is
  equal to $r \beta_{1,c+4}$. Hence we get the dimension formula of
  Theorem~\ref{refl3} from
  Theorem~\ref{thmeulerIF}  in this case. \\
  {\rm (iii)}. The lower bound of Theorem~\ref{thmeulerIF} is clearly better
  that the bound ${\rm ed} (\sF)$ provided $ \delta^2 > \delta^1 +
  \delta^3$.
\end{remark}

\begin{proof} By a general theorem of Laudal \cite{L1} on the dimension of the
  hull of any local deformation functor, we get that $ {_0\!\ext_R^1 }(F,F )-
  {_0\!\ext_R^2 }(F,F ) \le \dim O_F$ where $O_{F}$ is the hull of the
  deformation functor of the graded module $F$ (see Remark~\ref{graddefF}). To
  get the inequality of the theorem it suffices, by
  Proposition~\ref{propeulerIF}, to prove $ \dim O_{F} \le
  \dim_{(\sF)}{\rm M}_{\proj{3}}(c_1,c_2,c_3)$. Since we will use
  Theorem~\ref{localserre1} we replace $F$ by $F(v)$ for $v >>0$ to have the
  assumptions of Theorem~\ref{localserre1} satisfied. It is known the Hilbert
  scheme $\HH(d,g)$ contains a subscheme $H:=\HH(d,g)_{\gamma}$ which is the
  representing object of the subfunctor of flat families of curves with fixed
  postulation $\gamma$. For the local deformation functors at a curve $(X)$
  the latter corresponds precisely to the graded deformations of the
  homogeneous coordinate ring of $X$ (\cite{MDP1} and recall $\gamma(v) =
  h^0(\sI_X(v))$, $v \in \mathbb Z$, see also \cite{krao}). Hence we get
  \begin{equation} \label{ine1}
    \dim O_{H,(X)} = \dim_{(X)} \HH(d,g)_{\gamma} \le  \dim_{(X)}\HH(d,g).
\end{equation} 
By Theorem~\ref{localserre1}, $ \dim_{(\sF)}{\rm
  M}_{\proj{3}}(c_1,c_2,c_3)+h^0(\sF) = \dim_{(X)}\HH(d,g) +
h^0(\omega_X(-c_1+4))$. We {\it claim} that
\begin{equation} \label{ine2}
   \dim O_F+h^0(\sF) = \dim O_{H,(X)} + h^0(\omega_X(-c_1+4)).
\end{equation}
This is mostly explained in Remark~\ref{graddefF}. Indeed the natural
forgetful map $p_0: {\rm Def}_{F,H_{*}^0(\sigma)} \rightarrow {\rm Def}_{F}$
is smooth and has the same fiber as the forgetful map $ p: {\rm
  Def}_{\sF,\sigma} \rightarrow {\rm Def}_{\sF}$ in Theorem~\ref{localserre1}
by Remark~\ref{graddefF}. In the same way the corresponding graded variation
of $ q: {\rm Def}_{\sF,\sigma} \rightarrow {\rm Hilb}_{X/{\proj{3}}}$ is
smooth by $ _0\! \Ext^2(I_{X}(c_1),R) \simeq
\Ext^2(\sI_{X}(c_1),\sO_{\proj{3}})=0$ and its fiber coincides with that of
$q$, due to the isomorphism $ _0\! \Ext^1(I_{X}(c_1),R) \simeq
\Ext^1(\sI_{X}(c_1),\sO_{\proj{3}})$ (cf. \eqref{mspect} and \eqref{seq}) and
the arguments of Remark~\ref{graddefF}. This proves the claim and the
inequality of the theorem.

It remains to prove at ${_{-4}\!\hom_R}(F , M)$ is the correction number since
the reformulation as a kernel is trivial. Let $X$ be the generic curve of a
component of $\HH(d,g)$ which corresponds to $V$. Let $\gamma$ be the
postulation of $X$. Since there is a smooth open subscheme $U \ni (X)$ of
$\HH(d,g)$ of curves with postulation $\gamma$, we get $\HH(d,g)_{\gamma} \cap
U = \HH(d,g) \cap U$. Hence $\HH(d,g)_{\gamma}$ is smooth at $(X)$ and we have
equality in \eqref{ine1}. By Theorem~\ref{localserre1} and \eqref{ine2}, $
\dim O_{F} = \dim_{(\sF)}{\rm M}_{\proj{3}}(c_1,c_2,c_3)$ and $O_F$ is smooth.
Hence $\dim O_F= {_0\!\ext_R^1 }(F,F )$ and ${_0\!\ext_R^2 }(F,F )$ is the
correction number. Since we have ${_0\!\Ext_R^2 }(F,F )^{\vee} \simeq
{_{-4}\!\Hom_R}(F , M)$ by \eqref{mduality} and \eqref{mspect}, the proof is
complete.
\end{proof}
In \cite{krao} we proved a result (Lem.\! 2.2) similar to
Proposition~\ref{propeulerIF} for any curve $X$ with minimal resolution, 
\begin{equation} \label{resoluMX} 0 \rightarrow \bigoplus_i
  R(-i)^{\beta'_{3,i}} \rightarrow \bigoplus_i R(-i)^{\beta'_{2,i}}
  \rightarrow \bigoplus_i R(-i)^{\beta'_{1,i}} \rightarrow I_X \rightarrow 0\ ,
\end{equation}
namely that
\begin{equation} \label{lowerboundI}
  {_0\!\ext_R^1 }(I_X,I_X )- {_0\!\ext_R^2 }(I_X,I_X ) =  1 -
  \delta_I^0 = 4 d +  \delta_I^2 - \delta_I^1  \ 
\end{equation}
where $\delta_I^{j} = \sum _{i} (\beta'_{1,i} -\beta'_{2,i}+ \beta'_{3,i})
\cdot h^j(\sI_X(i))$ and $d= deg(X)$. Note that the difference of the $ {\ext
}$-numbers in \eqref{lowerboundI} is a lower bound for $ \dim
O_{\HH(d,g)_{\gamma},(X)}$ (\cite{krao}, proof of Thm.\! 2.6 (i)). As a
by-product of \eqref{ine1} and the proof above we get
\begin{theorem}\label{thmeulerIX}
  Let $X$ be a curve in $\proj{3}$. Then the dimension of $\dim_{(X)}\HH(d,g)$
  at $(X)$ satisfies
\begin{equation*}
  \dim_{(X)}\HH(d,g) \ge  1 -
  \delta_I^0 = 4 d +  \delta_I^2 - \delta_I^1.
\end{equation*}
Moreover if $X$ is a generic curve of a generically smooth component $V$ of $
\HH(d,g)$ and $M= H_{*}^1(\sI_X)$, then
$$\dim V = 4 d + \delta_I^2 
- \delta_I^1 +\ {_{-4}\!\hom_R}(I_X , M) $$ where
${_{-4}\!\Hom_R}(I_X , M)$ is the kernel of the map $ \bigoplus_i
H^1(\sI_X(i-4))^{\beta'_{1,i}} \rightarrow \bigoplus_i
H^1(\sI_X(i-4))^{\beta'_{2,i}} $ induced by \eqref{resoluMX}.
\end{theorem}

\begin{remark} Let $X$ be any curve in $\proj{3}$. \\
  {\rm (i)} If $M=0$, then $\delta_I^1=0$ and we can use Theorem 2.6 of
  \cite{krao} to see that the lower bound of Theorem~\ref{thmeulerIX} is equal
  to $ \dim_{(X)}\HH(d,g)$. 
  This coincides with \cite{E}.\\
  {\rm (ii)} If $\diam M=1$, $\dim M = r$ and $X$ is a generic unobstructed
  curve, then the lower bound is equal to $4 d + \delta_I^2 + r \beta'_{2,c}$
  because $ r \beta'_{1,c}=0$ for a generic curve by \cite{krao}, Cor.\! 4.4.
  Moreover in this case the ``correction'' number $ {_{-4}\!\hom_R}(I_X , M)$
  is equal to $r \beta'_{1,c+4}$. Hence we get $$\dim V = 4 d + \delta_I^2 +
  r( \beta'_{2,c} + \beta'_{1,c+4}). $$
  This coincides with the dimension formula of \cite{krao}, Thm.\! 3.4. \\
\end{remark}

\bigskip \bigskip


\begin{thebibliography}{10}
\setlength{\itemsep}{0pt}
\setlength{\parskip}{0pt}
{\small  

\bibitem{AO} V. Ancona, G. Ottaviani \newblock On singularities of ${\rm
    M}_{\proj{3}}(c_1,c_2)$. {\em Internat. J. Math. Vol.} 9 (1998), 407-419.
%
%
%
%

\bibitem{C} M-C. Chang.  \newblock A Filtered Bertini-type Theorem.  \newblock
  {\em  J. reine angew. Math.} 397 (1989), 214-219.

\bibitem{C1} M-C. Chang. \newblock Stable rank 2 reflexive sheaves on
  $\proj{}^3$ with large $c_3$. \newblock {\em J. reine angew.
    Math.} 343 (1983), 99--107.

\bibitem{C2} M-C. Chang. \newblock Stable rank 2 reflexive sheaves on
  $\proj{}^3$ with small $c_2$ and applications. \newblock {\em Trans. Amer.
    Math. Soc} 284 (1984), 57--84.


\bibitem{EF} Ph. Ellia, M. Fiorentini.  \newblock D\'{e}faut de postulation et
  singularit\'{e}s du Sch\'ema de Hilbert.  {\em Annali Univ. di Ferrara} 30
  (1984), 185-198.
%
%
\bibitem{ES} G. Ellingsrud, S.A. Strømme. \newblock Stable Rank-2 Vector
  Bundles with $c_1=0$ and $c_2=3$. \newblock {\em Math. Ann.} 255 (1981)
  123-135.

\bibitem{E} G. Ellingsrud. \newblock Sur le sch\'{e}ma de Hilbert des
  vari\'et\'es de codimension 2 dans $\proj{}^e$ a c\^{o}ne de
  Cohen-Macaulay, {\em Ann. Scient. \'{E}c. Norm. Sup.} {\bf 8} (1975),
  423-432.

%
%
\bibitem{Gu} S. Guffroy. \newblock Dimension des familles de courbes lisses
  sur une surface quartique normales de ${\proj{3}}$. \newblock To appear
  in {\em Proc. Amer. Math. Soc}.

\bibitem{GMir} P. Gurrola, R. M. Mir\'o-Roig. \newblock On the existence of
  generically smooth components for moduli spaces of rank 2 stable reflexive
  sheaves on $\proj{}^3$. \newblock {\em J. Pure Appl. Algebra}, 102 (1995),
  313--345.
     
\bibitem{F} G. Fløystad. \newblock Determining obstructions for space curves,
  with application to non-reduced components of the Hilbert scheme.  {\em J.
    reine angew. Math.} 439 (1993), 11-44.

\bibitem{SGA2} A. Grothendieck.  \newblock Cohomologie Locale des Faisceaux
  Cohérents et Théorèmes de Lefschetz Locaux et Globaux. Augmenté
  d'un exposé par M. Raynaud. $(SGA$ $2)$.  \newblock {\em Advanced
  Studies in Pure Mathematics} Vol. 2. North-Holland, Amsterdam (1968).

\bibitem{G} A. Grothendieck. \newblock Les sch\'{e}mas de Hilbert. \newblock
  {\em S\'{e}minaire Bourbaki}, exp. 221 (1960).
 
%
\bibitem{H1} R. Hartshorne.  \newblock Algebraic Geometry.  {\em Graduate
    Texts in Math.} Vol. 52 \newblock Springer--Verlag, New York, 1983.

\bibitem{H2} R. Hartshorne.  \newblock Stable vector bundles of rank 2 on
    $\proj{}^3$. {\em Math. Ann.} 238 (1978) 229-280.

  \bibitem{H3} R. Hartshorne. \newblock Stable reflexive sheaves. {\em Math.
      Ann.} 254 (1980) 121-176.

\bibitem{HS} R. Hartshorne, I. Sols. \newblock Stable rank 2 vector bundles
  on $\proj{}^3$ with $c_1 = -1$, $c_2 = 2$. {\em J. reine angew. Math.} 325
(1981), 145-152.

\bibitem{HL} D. Huybrechts, M. Lehn \newblock The Geometry of Moduli Spaces of
  Sheaves, Aspects of Math., E31, Vieweg (1997).
%
 %
%
%

\bibitem{IM} A. Iliev, D. Markushevich. \newblock Quartic 3-fold: Pfaffians,
  vector bundles, and half-canonical curves. {\em Michigan Math. J.} 47 no.2
  (2000), 385-394

\bibitem{K2} J. O. Kleppe. \newblock Deformations of reflexive sheaves of rank
  2 on $\proj{}^3$. Preprint March 1982, Univ. of Oslo.

\bibitem{K1} J. O. Kleppe.  \newblock Non-reduced components of the Hilbert
  scheme of smooth space curves, in \newblock {\em Proc. Rocca di Papa 1985}
  {\em Lectures Notes in Math.} Vol. 1266 \newblock Springer--Verlag (1987).
  

\bibitem{K3} J. O. Kleppe.  \newblock Liaison of families of subschemes in
  $\pp^n$, in \newblock ``Algebraic Curves and Projective Geometry,
    Proceedings (Trento, 1988),'' {\em Lectures Notes in Math.} Vol. 1389
  \newblock Springer--Verlag (1989).
 
\bibitem{K5} J. O. Kleppe.  \newblock Concerning the existence of nice
  components in the Hilbert scheme of curves in $\pp^n$ for $n=4$ and $5$,
 \newblock {\em J. reine angew. Math.} 475 (1996), 77-102.
 
\bibitem{krao} J. O. Kleppe. \newblock The Hilbert Scheme of Space Curves of
  small diameter. {\em Annales de l'institut Fourier} 56 no. 5 (2006),
  1297-1335.
   
  
  

  
 %
\bibitem{L1} A. Laudal.  \newblock Formal Moduli of Algebraic Structures.
  {\em Lectures Notes in Math.}, Vol. 754, \newblock Springer--Verlag, New
  York, 1979.

\bibitem{L} A. Laudal.  \newblock A generalized trisecant lemma, in \newblock
  Proceedings, Tromsø, 1977, {\em Lectures Notes in Math.}, Vol.
  687, \newblock Springer--Verlag, New York, 1978.
  
\bibitem{L2} A. Laudal.  \newblock Matrix Massey Products and Formal Moduli I,
  In: Algebra, Algebraic Topology and Their Interactions. {\em
    Lectures Notes in Math.}, Vol. 1183, 218-240 \newblock Springer--Verlag,
  New York, 1986.

\bibitem{MDP1} M. Martin-Deschamps, D. Perrin. \newblock Sur la classification
  des courbes gauches, {\em Asterisque} (1990), 184-185.
 
%
%
%
\bibitem{Mar1} M. Maruyama. \newblock Moduli of stable sheaves I, {\em J.
    Math.}, Kyoto Univ. 17 (1977), 91-166 .
  
\bibitem{Mar2} M. Maruyama. \newblock Moduli of stable sheaves II, {\em J.
    Math.}, Kyoto Univ. 18 (1978), 557-614.

\bibitem{MSS} J. Meseguer, I. Sols, S. A. Strømme, Compactification of a
  family of vector bundles on P3, In E. Balslev (Ed.): \newblock {\em 18th
    Scandinavian Congress of Mathematicians}, Proceedings, 1980, Birkhauser,
  Boston-Basel-Stuttgart 1981, 474-494.
  
\bibitem{Mig} J. Migliore.  \newblock Introduction to liaison theory and
  deficiency modules. {\em Progress in Math.}, Vol. 165, \newblock Birkhäuser
  Boston, Inc., Boston, MA, 1998.
%
%
%
\bibitem{Mir1} R. M. Mir\'o-Roig. \newblock Gaps in the Chern classes of rank 2
  stable reflexive sheaves, {\em Math. Ann.} 270 (1985), 317-323.

\bibitem{Mir2} R. M. Mir\'o-Roig.  \newblock Faisceaux reflexifs stables de
  rang 2 sur  $\proj{}^3$ non obstrues. {\em CRAS} 303 (1986), 711-713.
 
\bibitem{Mir3} R. M. Mir\'o-Roig. \newblock Some moduli spaces for rank 2
  stable reflexive sheaves on $\proj{}^3$. \newblock {\em Trans. Amer. Math.
    Soc} 299 (1987), 699--717.

%
\bibitem{Mir5} R. M. Mir\'o-Roig. \newblock Singular moduli spaces of stable
  vector bundles on $\proj{}^3$. \newblock {\em Pacific J. Math.} 172, no. 2
  (1996), 477--482.

\bibitem{MT} R. M. Mir\'o-Roig, G. Trautmann. \newblock The moduli scheme
  $M(0,2,4)$ over $\proj{}^3$. {\em Math. Z.} 216 (1994), 283--315.

%

\bibitem{Mu} D. Mumford.  \newblock Further pathologies in algebraic geometry,
 {\em Amer. J.  Math.}, 84, 1962, 642-648.
 
\bibitem{N} H. Nasu. \newblock Obstructions to deforming space curves and
  non-reduced components of the Hilbert scheme. {\em Publ. Res. Inst. Math.
    Sci} 42 (2006), 117--141.
  
%
\bibitem{R} A. P. Rao. \newblock Liaison Among Curves. {\em Invent. Math.} 50
  (1979), 205-217.

\bibitem{RVV} M. Roggero, P. Valabrega, M. Valenzano. \newblock Rank Two
  Bundles and Reflexive Sheaves on $\proj{}^3$ and Corresponding Curves: an
  Overview. in \newblock {\em Geometric and Combinatorial Aspects of
    Commutative Algebra (Messina 1999), LNPAM} 217, Marcel Dekker, 2001
  327-343.
  

\bibitem{Se} E. Sernesi. \newblock Un esempio di curva ostruita in
$\proj{3}$. {\em  Sem. di variabili Complesse, Bologna} (1981), 223-231.

\bibitem{St} S.A. Strømme. \newblock Ample Divisors on Fine Moduli Spaces on
the Projective Plane. {\em Math. Z.} 187 (1984), 405--423.

 %

\bibitem{Vak} R. Vakil. \newblock Murhpy's law in algebraic geometry: badly
  behaved deformation spaces. {\em Invent. Math.} 164 (2006), no. 3, 569-590.

\bibitem{Val} M. Valenzano. \newblock Rank 2 reflexive sheaves on a smooth
  threefold. \newblock {\em Rend. Sem. Mat. Univ. Pol. Torino} Vol 62, 3
  (2004).

\bibitem{Ver} P. Vermeire. \newblock Moduli of reflexive sheaves on smooth
  projective $3$-fold. J. Pure Appl. Algebra 211, No. 3, 622-632 (2007).

\bibitem{Vog} J.A. Vogelaar. \newblock Constructing vector bundles from
  codimension two subvarieties, PhD Thesis, Leiden (1978).

\bibitem{W1} C. Walter.  \newblock Some examples of obstructed curves in
  $\proj{3}$. \newblock {\em In: Complex Projective Geometry. London Math. Soc.
    Lecture Note Ser.} 179 (1992).
  
\bibitem{W2} C. Walter.  \newblock Horrocks theory and algebraic space
  curves.  Preprint (1990). 

\bibitem{Wei} C. Weibel. \newblock An Introduction to Homological Algebra.
  \newblock {\em Cambridge Studies in Advanced Mathematics 38, Cambridge
    University Press} (1994).
  \\

OSLO UNIVERSITY COLLEGE, FACULTY OF ENGINEERING, PB. 4 ST. OLAVS PLASS, N-0130
OSLO, NORWAY. 
 
 E-mail address: JanOddvar.Kleppe@iu.hio.no   \hspace{1.5cm} 
}
\end{thebibliography}
\end{document}